\documentclass[11pt]{article}
\usepackage{amsmath,amsthm,amssymb,color,epic,eepic, amsfonts,graphicx,
cite
}
\renewcommand{\baselinestretch}{1.2}

\def\baselinestretch{1.4}
\newlength{\minitwocolumn}
\newcommand{\Z}{{\Bbb Z}} 
\newcommand{\C}{{\Bbb C}} 
\newcommand{\N}{{\Bbb N}} 
\newcommand{\FF}{{\Bbb F}} 
\newcommand{\PP}{{\Bbb P}} 
\newcommand{\F}{{\mathcal F}}
\newcommand{\A}{{\mathbb A}}

\newcommand{\cD}{{\mathcal D}}
\newcommand{\cA}{{\mathcal A}}
\newcommand{\cU}{{\mathcal U}}
\newcommand{\cZ}{{\mathcal Z}}
\newcommand{\cB}{{\mathcal B}}

\newcommand{\cH}{{\mathcal H}}
\newcommand{\cN}{{\mathcal N}}

\newcommand{\cR}{{\mathcal R}}

\newcommand{\cM}{{\mathcal M}}

\newcommand{\cQ}{{\mathcal Q}}

\newcommand{\cV}{{\mathcal V}}
\newcommand{\cW}{{\mathcal W}}

\renewcommand{\H}{{\mathcal H}}
\newcommand{\la}{\lambda}

\newcommand{\si}{{\sigma}}

\newcommand{\al}{\alpha}

\newcommand{\vep}{\varepsilon}

\newcommand{\s}{{\sigma}}

\newcommand{\hpi}{{\pi}}
\newcommand{\hf}{\widehat{f}}

\newcommand{\hV}{{\cV}}

\newcommand{\nn}{{\nonumber}}
\newcommand{\bea}{\begin{eqnarray}}
\newcommand{\ena}{\end{eqnarray}}
\newcommand{\be}{\begin{eqnarray*}}
\newcommand{\en}{\end{eqnarray*}}

\newcommand{\beit}{\begin{itemize}}
\newcommand{\enit}{\end{itemize}}

\newcommand{\lb}[1]{\label{#1}}

\newcommand{\ds}[1]{{\displaystyle #1 }}

\newcommand{\End}{\mathrm{ End}}
\newcommand{\rank}{\mathrm{ rank}}

\newcommand{\Ind}{\mathrm{ Ind}}
\newcommand{\id}{\mathrm{ id}}

\newcommand{\wt}{\mathrm{ wt}}

\newcommand{\E}{\mathrm{E}}

\newcommand{\gC}{\mathfrak{C}}
\newcommand{\gD}{\mathfrak{D}}
\newcommand{\gN}{\mathfrak{N}}
\newcommand{\gH}{\mathfrak{H}}

\newcommand{\gS}{\mathfrak S}

\newcommand{\ket}[1]{{| #1 \rangle}}      

\def\infq4p#1{{(#1;q^4,p)_\infty}}

\newcommand{\tot}{\, \widetilde{\otimes}\, }

\newcommand{\mmatrix}[1]{\begin{matrix} #1 \end{matrix}}

\font\teneufm=eufm10
\font\seveneufm=eufm7
\font\fiveeufm=eufm5
\newfam\eufmfam
\textfont\eufmfam=\teneufm
\scriptfont\eufmfam=\seveneufm
\scriptscriptfont\eufmfam=\fiveeufm
\def\frak#1{{\fam\eufmfam\relax#1}}
\let\goth\frak
\newcommand{\slth}{\widehat{\goth{sl}}_2}
\newcommand{\slt}{\goth{sl}_2}

\newcommand{\slnh}{\widehat{\goth{sl}}_N}
\newcommand{\sln}{\goth{sl}_N}
\newcommand{\g}{\goth{g}}
\newcommand{\gt}{\goth{g}_{tor}}

\newcommand{\Bqla}{{{\mathcal B}_{q,\lambda}}}

\newcommand{\gl}{{\goth{gl}}}
\newcommand{\gln}{{\goth{gl}_N}}

\newcommand{\glnh}{\widehat{\goth{gl}}_N}

\newcommand{\h}{\goth{h}}
\newcommand{\gh}{\widehat{\goth{g}}}

\font\fourteeneufm=eufm10 scaled\magstep2    
   
\newcommand{\gbig}{\mbox{\fourteeneufm g}} 
\makeatletter
\@addtoreset{equation}{section}
\makeatother

\newtheorem{thm}{Theorem}[section]
\newtheorem{prop}[thm]{Proposition}
\newtheorem{lem}[thm]{Lemma}
\newtheorem{cor}[thm]{Corollary}
\newtheorem{conj}[thm]{Conjecture}

\newtheorem{df}{Definition}[section]
\newtheorem{dfn}[thm]{Definition}

\newcommand{\ip}{{i'}} 
\newcommand{\ipp}{{i''}} 
\newcommand{\spp}{{s'}} 

\begin{document}

\vspace{-1cm}
\begin{center}
{\bf\Large  Elliptic Quantum Toroidal Algebras, $Z$-algebra Structure\\ and  Representations
\\[7mm] }
{\large  Hitoshi Konno${}^{\dagger}$ and Kazuyuki Oshima${}^{\star}$ }\\[6mm]
${}^\dagger${\it  Department of Mathematics, Tokyo University of Marine Science and 
Technology, \\Etchujima, Koto, Tokyo 135-8533, Japan\\
       hkonno0@kaiyodai.ac.jp}\\
${}^\star${\it Center for General Education, Aichi Institute of Technology, \\ Yakusa-cho, Toyota
470-0392, Japan }\\
oshima@aitech.ac.jp
\\[7mm]
\end{center}

\begin{abstract}
We introduce a new elliptic quantum toroidal algebra $U_{q,\kappa,p}(\g_{tor})$ associated with an  arbitrary toroidal algebra $\g_{tor}$. 
We show that $U_{q,\kappa,p}(\g_{tor})$ contains 
two elliptic quantum algebras associated with a corresponding affine Lie algebra  $\gh$ as subalgebras. 
They are analogue of the horizontal and the vertical subalgebras in the quantum toroidal algebra $U_{q,\kappa}(\g_{tor})$. 
A Hopf algebroid structure is introduced as a co-algebra structure of  $U_{q,\kappa,p}(\g_{tor})$ using the Drinfeld comultiplication. 
We also investigate the $Z$-algebra structure of $U_{q,\kappa,p}(\g_{tor})$ and show that the 
$Z$-algebra governs the irreducibility of the level $(k(\not=0),l)$-infinite dimensional $U_{q,\kappa,p}(\g_{tor})$-modules in the same way as in the elliptic quantum group $U_{q,p}(\gh)$. As an example,  we construct the level $(1,l)$ irreducible representation of $U_{q,\kappa,p}(\g_{tor})$ for the simply laced $\g_{tor}$. 
We also construct the level  $(0,1)$  representation of $U_{q,\kappa,p}(\gl_{N,tor})$ and discuss a conjecture on its geometric interpretation as an action   on the torus equivariant elliptic cohomology of the affine $A_{N-1}$ quiver variety. 

%
%
\end{abstract}

\section{Introduction}
Quantum toroidal algebras were introduced in \cite{GKV}. Let $\g_{tor}$ denote a toroidal algebra associated with a complex semisimple Lie algebra $\g$, i.e. 
$\g_{tor}$ is a canonical two-dimensional central extension of the double loop Lie algebra $\C^{\times}\times \C^{\times} \to \g$\cite{GKV}. We denote by 
$U_{q,\kappa}(\g_{tor})$ the quantum toroidal algebra associated with $\g_{tor}$. Let $\gh$ denote an untwisted affine Lie algebra associated with $\g$, i.e. 
$\gh$ is a one-dimensional canonical central extension of the loop Lie algebra $\C^{\times}\to \g$.   The quantum toroidal algebra $U_{q,\kappa}(\g_{tor})$ was formulated as an extension of 
the quantum affine algebra $U_q({\gh})$  in  the Drinfeld  realization\cite{Drinfeld} by replacing  the finite type Cartan matrix in the latter  with the affine type generalized Cartan matrix.  For the case $\g=\gl_N$ 
one can introduce an extra parameter $\kappa$\cite{GKV} due to a cyclic property of the affine Dynkin diagram of $A_{N-1}$ type.    
Representation theory of $U_{q,\kappa}(\g_{tor})$ has been developed by many works\cite{FJMM2, FJMM4,FJM, Na01,VV98, VV99,Nagao, Miki1, Miki2, Miki3, STU, Saito,H}.  Among others, we would like to mention the following results in the case $\g=\gln$ with $N\geq 3$: a toroidal analogue of the Schur duality\cite{VV96} between $U_{q,\kappa}(\gl_{N,tor})$ and the double affine Hecke algebra\cite{Ch}, $q$-Fock space representations\cite{VV98, STU,FJMM2}, a family of quasi-finite representations including the Macmahon modules\cite{FJMM2}. 
Their application to a calculation of the instanton partition function on the ALE space was also discussed in \cite{AKMMSZ}. 
See also a development of representations in terms of the shuffle algebra\cite{Negut,Tsym}.  

On the other hand, the case $\g=\gl_1$ the quantum toroidal algebra $U_{q,t}(\g_{1,tor})$  was introduced in \cite{Miki4} as a $q,t$-deformation of the $W_{1+\infty}$ algebra, and 
various representations have been studied \cite{AFS,Miki4, FT, FHHSY, FFJMM, FJMM1, FJMM3, BFM}.  It is remarkable  that  $U_{q,t}(\gl_{1,tor})$ is  isomorphic to the elliptic Hall algebra\cite{SV1,SV2,Sc1,FFJMM}.  
A deep connection of representations of  $U_{q,t}(\gl_{1,tor})$ to the Macdonald theory was clarified\cite{Miki4,FT, FHHSY, FFJMM,SV1,FOS19}. 
Representations of  $U_{q,t}(\gl_{1,tor})$ were also applied  
to the 5d and 6d lifts of the 4d ${\cN=2}$ SUSY gauge theories known as 
the linear quiver gauge theories, 
such as a calculation of the instanton partition functions\cite{Nekrasov04} and a study of  
the higher dimensional analogue of the Alday-Gaiotto-Tachikawa (AGT) correspondence\cite{AGT}. 
See for example \cite{AFS,MMZ,Zenke,AKMMMMOZ,BFMZZ,BFHMZ,Nieri}. 
There are some different approaches based on the elliptic Hall algebra \cite{SV3} and the shuffle algebra\cite{NegutAGT
}. For the affine Yangian case see \cite{MO}. 

In the previous paper\cite{KO}, we introduced the elliptic quantum toroidal algebra 
$U_{q,t,p}(\gl_{1,tor})$ as an elliptic analogue of $U_{q,t}(\gl_{1,tor})$. Constructing  representations
including the level (0,0) representation realized in terms of the elliptic Ruijsenaars difference operator, we discussed a deep connection of  $U_{q,t,p}(\gl_{1,tor})$ with an expected elliptic analogue of the Macdonald symmetric functions. 
We also showed that $U_{q,t,p}(\gl_{1,tor})$ realizes 
 the affine quiver $W$ algebra $W_{q,t}(\Gamma(\widehat{A}_0))$ proposed in  \cite{KimPes}. Namely, the elliptic currents and a composition of the vertex operators of $U_{q,t,p}(\gl_{1,tor})$ give a realization of the screening currents 
and the generating currents of $W_{q,t}(\Gamma(\widehat{A}_0))$, respectively.   
It is also remarkable  that such realization provides a   
  relevant scheme to the  instanton calculus of  
 the 5d and 6d lifts of the  4d $\cN=2^*$ SUSY gauge theories, i.e.  the 
  $\cN=2$ SUSY gauge theories coupled with the adjoint matter\cite{Nekrasov04,Nekrasov}, known as  
 the Jordan quiver gauge theories\footnote{It is also called the ADHM quiver gauge theories.}.  
Hence we obtained 
a new AGT type correspondence.  
Here one of the essential feature is that $U_{q,t,p}(\gl_{1,tor})$  possesses  the  four parameters $q, t, p, p^*$ with one constraint  $p/p^*=t/q$ in  the level  $(1,N)$ representation, which play the role of  the $SU(4)$ $\Omega$ deformation parameters\cite{Nekrasov}.

A connection of the elliptic quantum algebras to the deformed $W$ algebras is  quite natural. For example, the elliptic quantum group $U_{q,p}(\gh)$ can be regarded as 
a $q$-deformation of the Feigin-Fuchs construction\cite{FeiFuch}, i.e. 
a systematic procedure of  constructing $W$ algebra as the Goddard-Kent-Olive (GKO) coset theory $\gh_{r-h^\vee-k}\oplus \gh_k\supset \gh_{r-h^\vee}$\cite{GKO,LF,BNY,KMQ,DJKMO,CrRav,Rav} by considering a $r$-deformation of $\gh_k$. 
 Here  $k$ denotes a level of representation of 
$\gh$, and  a $q$-deformation of $r$ becomes the {\it elliptic} nome $p=q^{2r}$ in $U_{q,p}(\gh)$.  Then the ordinary 
deformation parameters $q, t$ are identified with $p, p^*=pq^{-2k}$ in $U_{q,p}(\gh)$. See \cite{K98,KonnoBook}\footnote{It is also interesting to remark that  the elliptic quantum group $U_{q,p}(\glnh)$ was essentially realized in the tensor product of the quantum toroidal algebra $U_{q,t}(\gl_{1,tor})$\cite{BFM,FOS19,FOS20}. }. 
The importance of considering the elliptic quantum group $U_{q,p}(\gh)$ lies in  the fact that it provides an algebraic structure  i.e. a co-algebra structure,  which enables us to define and construct the vertex operators of the deformed $W$-algebras 
as intertwining operators of the $U_{q,p}(\gh)$-modules\cite{Konno08,KonnoBook}.

The aim of this paper is to introduce a new elliptic quantum toroidal algebra  $U_{q,\kappa,p}(\g_{tor})$ as a 
$p$-deformation of the quantum toroidal algebra $U_{q,\kappa}(\g_{tor})$. 
We formulate $U_{q,\kappa,p}(\g_{tor})$  by generators and relations in the same scheme as the elliptic quantum group $U_{q,p}(\gh)$
\cite{K98,JKOS,FKO,Konno18,KonnoBook}. 
In particular, we clarify the  $Z$-algebra structure of $U_{q,\kappa,p}(\g_{tor})$, which has not yet investigated even in the trigonometric case $U_{q,\kappa}(\g_{tor})$ systematically. We show that 
the  $Z$-algebras of $U_{q,\kappa,p}(\g_{tor})$ and  $U_{q,\kappa}(\g_{tor})$ are the same except for that  the former  depends on 
the dynamical parameters. 
This feature should be essential to $U_{q,\kappa,p}(\g_{tor})$ having a characterization as a $q$-deformation of  the Feigin-Fuchs construction i.e. a $q$-deformed $W$ algebra. See the case $U_{q,p}(\gh)$\cite{FKO}.  
We also construct  important representations such as the level $(1,l)$ ones for simply laced $\g_{tor}$ and the $(0,1)$ representation for $\gl_{N,tor}$ with $N\geq 3$.  
 In separate publications, we will show that $U_{q,\kappa,p}(\g_{tor})$ provides a realization of a $q$-deformation of the quiver $W$-algebra $W_{p,p^*}(\Gamma)$ associated with the affine Dynkin diagram $\Gamma$\cite{KimPes}.

This paper is organized as follows. 
In Sec.2, after reviewing basic facts in  the quantum toroidal algebra $U_{q,\kappa}(\g_{tor})$, we give a definition of the elliptic quantum toroidal algebra $U_{q,\kappa,p}(\g_{tor})$. We show that $U_{q,\kappa,p}(\g_{tor})$ contains 
two elliptic quantum algebras isomorphic to $U_{q,p}(\gh)$ as subalgebras. 
These are analogue of the horizontal and the vertical subalgebras in $U_{q,\kappa}(\g_{tor})$\cite{GKV}.  In Sec.3 we introduce the Hopf algebroid structure as a co-algebra structure of  $U_{q,\kappa,p}(\g_{tor})$ using 
the Drinfeld comultiplication. In Sec.4, we define the dynamical representation of 
$U_{q,\kappa,p}(\g_{tor})$. In Sec.5, we discuss 
the $Z$-algebra structure. We show that the $Z$-algebra governs the irreducibility of the infinite dimensional $U_{q,\kappa,p}(\g_{tor})$-modules in the same way as in $U_{q,p}(\gh)$\cite{FKO}. As an example,  we give the level $(1,l)$ irreducible representation of $U_{q,\kappa,p}(\g_{tor})$ 
for the simply laced $\g_{tor}$. 
In Sec.6 we construct the level  $(0,1)$  representation of $U_{q,\kappa,p}(\gl_{N,tor})$ and discuss a conjecture on its geometric interpretation as an action of $U_{q,\kappa,p}(\gl_{N,tor})$  on the torus equivariant elliptic cohomology of the affine $A_{N-1}$ quiver variety. 
Appendix A is a list of formulas, which is used to discuss the $Z$-algebra structure in Sec.\ref{Zalgstr}. Appendix B is a direct proof of  the level $(0,1)$ representation of $U_{q,\kappa,p}(\gl_{N,tor})$.

A part of the results in Sec.2, 3, 4, 6 have been presented by H.K. at several workshops \cite{K19}.

\section{Elliptic Quantum Toroidal Algebras $U_{q,\kappa,p}(\gbig_{tor})$ }\lb{app:1}
After reviewing the quantum toroidal algebra $U_{q,\kappa}(\g_{tor})$\cite{GKV,SaitoLec}, we introduce  the elliptic quantum toroidal algebra  $U_{q,\kappa,p}(\g_{tor})$. Through this paper $q, p, \kappa \in \C^{\times}$ are generic complex numbers. 

\subsection{Quantum toroidal algebras $U_{q,\kappa}(\g_{tor})$}
Let $\g$ be a simple Lie algebra over $\C$ and $\gh$  the corresponding untwisted affine Lie algebra with the generalized Cartan 
matrix $A=(a_{ij})_{i,j \in I}$, $I=\{0,1,\dots,N\}$, \ $\rank A=N$. 
We denote by $B=(b_{ij})_{i,j \in I}, \,b_{ij}=d_i a_{ij}$ the symmetrization of $A$. 
We set $q_i=q^{d_i}$. We also use the notations for $n \in \mathbb{Z}$,
\begin{eqnarray}
&&[n]=\frac{q^n-q^{-n}}{q-q^{-1}},\\
&&[n]_i=\frac{q_i^n-q_i^{-n}}{q_i-q_i^{-1}}, \quad [n]_i !=[n]_i \cdots [1]_i, \quad \left[ n \atop k \right]_i=\frac{[n]_i !}{[n-k]_i ! [k]_i !}.
\end{eqnarray}
We fix a realization $(\h,\Pi,\Pi^\vee)$ of $A$, that is, $\h$ is a $N+2$-dimensional $\mathbb{C}$-vector space, 
$\Pi=\{\alpha_0,\alpha_1, \dots, \alpha_N\} \subset \h^*$ a set of simple roots, and $\Pi^\vee=\{h_0,h_1,\dots,h_N\} \subset \h$ a set of simple coroots satisfying 
$\langle\alpha_j,h_i\rangle=a_{ij}\ (i,j\in I)$ for a canonical pairing $\langle\, ,\rangle : \h^*\times \h \to \C$.  We also set $\cQ=\sum_{i\in I}\Z \al_i$. We take $\{h_1,\dots,h_N, c, d\}$  as the basis of $\h$ and $\{ \bar{\Lambda}_1,\cdots, \bar{\Lambda}_N, \Lambda_0, \delta \}$ the dual basis satisfying
\bea
&&
\langle\delta,d\rangle=1=\langle\Lambda_0,c\rangle,\quad 
\langle\bar{\Lambda}_i,h_j\rangle=\delta_{i,j}\lb{pairinghhs}, 
\ena
with the other pairings being 0. 
Denote also by $(a^\vee_0,a^\vee_1,\cdots, a^\vee_N)$ the co-label of the affine Dynkin diagram\cite{Kac}.

Let $M=(m_{ij})_{i,j \in I}$ be a matrix defined by 
\bea
&&m_{ij}=\left\{\mmatrix{\delta_{i,j+1}-\delta_{i+1,j}
&\mbox{for}\ \g=A_N \cr
0 &\mbox{otherwise} \cr}\right. .
\ena

\begin{df}\cite{GKV}\lb{defQTA}
Let $\kappa \in \mathbb{C}^\times$ be an additional parameter. 
The quantum toroidal algebra $U_{q,\kappa}(\g_{tor})$ is a unital  $\mathbb{C}$-algebra generated by 
\be
a^\vee_{i,l}, \quad X^+_{i,n}, \quad X^-_{i,n}, \quad k_i^{\pm1}, \quad q^{\pm c/2}, \quad q^d\quad 
 (i\in I, \quad l \in \mathbb{Z}\backslash\{0\}, \quad n \in \mathbb{Z}). 
\en
Set the generating functions called the Drinfeld currents as 
\be
&&X^+_i(z)=\sum_{n\in \Z}X^+_{i,n} z^{-n},\qquad X^-_i(z)=\sum_{n\in \Z}X^-_{i,n} z^{-n},\\
&&\Phi^\pm_i(z)=
k_i^{\pm1} \exp\left( \pm(q_i-q_i^{-1}) \sum_{n>0} a^\vee_{i,\pm n}z^{\mp n}\right).
\en
Then the  defining relations are given in terms of the  Drinfeld currents as follows. 
\bea
&& q^{\pm c/2},\quad \prod_{i\in I}k^{\pm a^\vee_i}_i\  :\hbox{ central },\label{qta1}\\
&& k_ik_{i}^{-1}=1=k_{i}^{-1}k_i,\quad  [k_i,k_j]=0=[k_i,a^\vee_{j,k}], \label{qta2}\\ 
&&[q^d, k^{\pm 1}_i]=0,\quad [q^d, a^\vee_{i,l}]=q^la^\vee_{i,l},\quad [q^d, X^\pm_{i}(z)]=X^\pm_i(q^{-1}z),\\
&& k_i X^+_j(z)k_{i}^{-1}=q^{\langle \al_j,h_i \rangle }X^+_j(z), \quad k_i X^-_j(z)k_{i}^{-1}=q^{-\langle\al_j,h_i\rangle }X^-_j(z),\label{qta3}\\
&&[a^\vee_{i,l},a^\vee_{j,m}]=\delta_{l+m,0}\frac{[a_{ij} l]_i}{l}\frac{q^{lc}-q^{-lc}}{q_j-q_j^{-1}}\kappa^{-lm_{ij}}, \label{qta4}\\
&&
[a_{i,l},X^+_j(z)]=\frac{[a_{ij} l]_i}{l}q^{-|l|c/2}\kappa^{-lm_{ij}}z^l X^+_j(z),\\
&&
[a_{i,l},X^-_j(z)]=-\frac{[a_{ij} l]_i}{l} q^{|l|c/2} \kappa^{-lm_{ij}}z^l X^-_j(z),
\\
&&(\kappa^{m_{ij}}z-q^{b_{ij}}w)
X^+_i(z)X^+_j(w)= (\kappa^{m_{ij}}q^{ b_{ij}}z-w) X^+_j(w)X^+_i(z),
\\
&&(\kappa^{m_{ij}}z-q^{-b_{ij}}w)
X^-_i(z)X^-_j(w)= (\kappa^{m_{ij}}q^{- b_{ij}}z-w) X^-_j(w)X^-_i(z),
\\
&&[X^+_i(z),X^-_j(w)]=\frac{\delta_{i,j}}{q_i-q_i^{-1}}
\left(\delta\bigl(q^c\frac{w}{z}\bigr)\Phi^+_i(q^{\frac{c}{2}}w)
-\delta\bigl(q^c\frac{z}{w}\bigr)\Phi^-_i(q^{\frac{c}{2}} z)
\right),
\ena
\bea
&&\sum_{\sigma\in \gS_a}\sum_{r=0}^{a} (-1)^r
\left[\mmatrix{a \cr
 r\cr}\right]_i
X^+_i(z_{\sigma(1)})\cdots X^+_i(z_{\sigma(r)})
X^+_{j}(w) X^+_i(z_{\sigma(r+1)})\cdots X^+_i(z_{\sigma(a)})=0,
\nonumber
\\
&&\\
&&\sum_{\sigma\in \gS_a}\sum_{r=0}^{a} (-1)^r
\left[\mmatrix{a \cr
 r\cr}\right]_i
 X^-_i(z_{\sigma(1)})\cdots X^-_i(z_{\sigma(r)})
X^-_{j}(w) X^-_i(z_{\sigma(r+1)})\cdots X^-_i(z_{\sigma(a)})=0,
\nonumber\\
&&
\ena
where we set $a=1-a_{ij}$ and $\delta(z)=\sum_{n\in\Z}z^n$.

\end{df}

\begin{prop}\cite{GKV}\lb{propGKV}
The quantum toroidal algebra ${U}_{q,\kappa}(\g_{tor})$ has the quantum affine algebra  
${U}_{q}(\gh)$ generated by $a^\vee_{i,l}, k^{\pm1}_i, X^+_{i,n}, X^-_{i,n}, q^{\pm {c}/{2}}$ $(i\in I\backslash \{0\},  
l\in \Z\backslash\{0\}, n\in \Z )$ as a subalgebra. We call it the horizontal subalgebra.
${U}_{q,\kappa}(\g_{tor})$ has another subalgebra  
${U}_{q}(\gh)$ generated by $X^+_{i,0}, X^-_{i,0},  k^{\pm1}_i$ $(i\in I)$, which is again a quantum affine algebra.  We call it the vertical subalgebra. 
\end{prop}

\subsection{The elliptic quantum toroidal algebra $U_{q,\kappa,p}(\g_{tor})$}
Let $H_P$ be  a $\C$-vector space spanned by $ P_0, P_1,\cdots,P_N$
 and $H^{*}_P$ be its dual space spanned by $Q_0, Q_1,\cdots,Q_N$ with a 
 pairing $\langle Q_i,P_j\rangle =a_{ij}$.   For $\al=\sum_{i\in I}c_i\al_i\in \h^*$, we set $Q_\al=\sum_{i\in I}c_iQ_i$. In particular $Q_{\al_i}=Q_i$. Define also an analogue of the root lattice $\cR_Q=\sum_{i\in I}\Z Q_i$. 
 
Let us consider $H:=\h\oplus H_P$ and  $H^*:=\h^*\oplus H^*_P$ 
 with a pairing $\langle \ ,\ \rangle : H^*\times H \to \C$ by extending those on $\h^*\times \h$ and $H^*_P\times H_P$ with $\langle \h^* ,H_P \rangle=0=\langle H^*_P ,\h \rangle$. 
We denote by $\FF=\cM_{H^*}$ the field of meromorphic functions on $H^*$. We regard 
a meromorphic function $g(h,P)$ of $h, P\in H$ as an element in $\cM_{H^*}$ by 
$g(h,P)(\mu)=g(\langle \mu,h\rangle, \langle \mu,P\rangle)$ for $\mu\in H^*$.

\begin{dfn}\lb{defEQTA}
The elliptic quantum toroidal algebra $U_{q,\kappa,p}(\g_{tor})$ is a topological algebra over  $\FF[[p]]$  
generated by 
\be
\al^\vee_{i,l}, \quad x^{\pm}_{i,n}, \quad K_i^{\pm}, \quad q^{\pm c/2}, \quad q^{d}\quad 
 (i\in I, \quad l \in \mathbb{Z}\backslash\{0\}, \quad n \in \mathbb{Z}). 
\en
The defining  relations can be  written in terms of  the generating functions called the elliptic currents  given by
\be
&&x^{\pm}_i(z)=\sum_{n\in \Z}x^{\pm}_{i,n} z^{-n},\\
&&\phi^+_i(q^{\frac{c}{2}}z)
=K^+_i \exp\left( -(q_i-q_i^{-1}) \sum_{n>0} \frac{p^n \al^\vee_{i, -n}}{1-p^n}z^{ n}\right)
 \exp\left( (q_i-q_i^{-1}) \sum_{n>0} \frac{\al^\vee_{i, n}}{1-p^n}z^{ -n}\right),\\
&&\phi^-_i(q^{-\frac{c}{2}}z)
=K^-_i \exp\left( -(q_i-q_i^{-1}) \sum_{n>0} \frac{ \al^\vee_{i, -n}}{1-p^n}z^{ n}\right)
 \exp\left( (q_i-q_i^{-1}) \sum_{n>0} \frac{p^n \al^\vee_{i, n}}{1-p^n}z^{ -n}\right).
\en
We also set $p^*=pq^{-2c}$. 
The defining relations are given as follows. 
For $g(h,P) \in \cM_{H^*}$, 
\bea
&& q^{\pm c/2}, \quad \prod_{i\in I}(K_i^{\pm})^{ a^\vee_i}\  :\hbox{ central },\\
&&
g(h,{P})x^+_j(z)=x^+_j(z)g(h+\langle \al_j,h\rangle, P-\langle Q_{j},P\rangle),\quad
g(h,{P})x^-_j(z)=x^-_j(z)g(h-\langle \al_j,h\rangle, P),\nn\\
&&\lb{gf}\\
&&[g(P), \al^\vee_{i,m}]=
0,\qquad [g(h, P), q^d]=0,\quad 
g(h,{P})K^{\pm}_j=K^{\pm}_jg(h, P-\langle Q_{j},P\rangle),
\lb{gKpm}
\\
&& [q^d, \al^\vee_{j,m}]=q^m\al^\vee_{j,m},\quad [q^d, x^\pm_j(z)]=x^\pm_j(q^{-1}z)
\lb{dedf}\\
&&
[K^\pm_i,K^\pm_j]=[K^\pm_i,K^\mp_j]=0=[K^\pm_i,\al^\vee_{j,m}], \label{qta2}\\ 
&& K^\pm_i x^+_j(z)(K^\pm_{i})^{-1}=q^{\pm\langle \al_j,h_i \rangle }x^+_j(z), \quad K^\pm_i x^-_j(z)(K^\pm_{i})^{-1}=q^{\mp\langle\al_j,h_i\rangle }x^-_j(z),\label{qta3}\\
&&[\al^\vee_{i,l},\al^\vee_{j,m}]=\delta_{l+m,0}\frac{[a_{ij} l]_i}{l}\frac{q^{lc}-q^{-lc}}{q_j-q_j^{-1}}\frac{1-p^l}{1-p^{*l}}
\kappa^{-lm_{ij}}q^{-lc},\\
&&
[\al^\vee_{i,l},x^+_j(z)]=\frac{[a_{ij} l]_i}{l}\frac{1-p^l}{1-p^{*l}}q^{-lc}\kappa^{-lm_{ij}}z^l x^+_j(z),\\
&&
[\al^\vee_{i,l},x^-_j(z)]=-\frac{[a_{ij} l]_i}{l}  \kappa^{-lm_{ij}}z^l x^-_j(z),
\ena
\bea
&&(\kappa^{m_{ij}} z-q^{b_{ij}}w)g_{ij}(\kappa^{-m_{ij}}w/z;p^*)
x^+_i(z)x^+_j(w)
=(\kappa^{m_{ij}}q^{b_{ij}}z-w) g_{ij}(\kappa^{m_{ij}}z/w;p^*)
x^+_j(w)x^+_i(z),\nn\\
&&\lb{xpxpgen}\\
&&(\kappa^{m_{ij}} z-q^{-b_{ij}}w)g_{ij}(\kappa^{-m_{ij}}w/z;p)^{-1}
 x^-_i(z)x^-_j(w)
=(\kappa^{m_{ij}}q^{-b_{ij}}z-w) g_{ij}(\kappa^{m_{ij}}z/w;p)^{-1}
x^-_j(w)x^-_i(z),\nn\\
&&\lb{xmxmgen}\\
&&[x^+_i(z),x^-_j(w)]=\frac{\delta_{i,j}}{q_i-q_i^{-1}}
\left(\delta\bigl(q^c \frac{w}{z}\bigr)\phi^{+}_{i}(q^{\frac{c}{2}} w)
-\delta\bigl(q^{-c} \frac{w}{z}\bigr)\phi^{-}_{i}(q^{\frac{c}{2}}z)
\right),\lb{xpxm}\\
&&\sum_{\sigma\in \gS_a}\prod_{1\leq k<m\leq a
}g_{ii}(z_{\sigma(m)}/z_{\sigma(k)};p^*)
\nonumber\\
&&\qquad\quad\times
\sum_{r=0}^a (-1)^r
\left[\mmatrix{a\cr r\cr}\right]_i \prod_{1\leq s\leq r}
g_{ij}(\kappa^{-m_{ij}}w/z_{\sigma(s)};p^*)
\prod_{r+1\leq s\leq a} 
g_{ij}(\kappa^{m_{ij}}z_{\sigma(s)}/w;p^*)
\nn \\
&&\qquad\quad\times x^+_{i}(z_{\sigma(1)})\cdots x^+_{i}(z_{\sigma(r)})
x^+_{j}(w) x^+_{i}(z_{\sigma(r+l)})\cdots x^+_{i}(z_{\sigma(a)})=0,\lb{Serrexp}\\
&&\sum_{\sigma\in \gS_a}\prod_{1\leq k<m\leq a
}g_{ii}(z_{\sigma(m)}/z_{\sigma(k)};p)^{-1}
\nonumber\\
&&\qquad\quad\times
\sum_{r=0}^a (-1)^r
\left[\mmatrix{a\cr r\cr}\right]_i \prod_{1\leq s\leq r}
g_{ij}(\kappa^{-m_{ij}}w/z_{\sigma(s)};p)^{-1}
\prod_{r+1\leq s\leq a} 
g_{ij}(\kappa^{m_{ij}}z_{\sigma(s)}/w;p)^{-1}
\nn\\
&&\qquad\quad\times x^-_{i}(z_{\sigma(1)})\cdots x^-_{i}(z_{\sigma(r)})
x^-_{j}(w) x^-_{i}(z_{\sigma(r+l)})\cdots x^-_{i}(z_{\sigma(a)})=0.
\ena
Here $a=1-a_{ij}$ as before, and we set 
\bea
&&g_{ij}(z;s)=\exp\left(-\sum_{m>0}\frac{1}{m}\frac{q^{b_{ij}m}-q^{-b_{ij}m}}{1-s^m}(sz)^m\right)\ \in \ \C[[s]][[z]].\lb{deg:g}
\ena
We also denote by ${U}'_{q,\kappa,p}(\g)$  the same elliptic quantum toroidal algebra 
 ${U}_{q,\kappa,p}(\g)$ without the element $q^d$. 
\end{dfn}
We treat these relations as formal Laurent series in the argument of the elliptic currents i.e. $z, w, z_1, \cdots, z_a$. 
Then all the coefficients in $z,w$ etc. are well defined in the $p$-adic topology\cite{Konno18}.

It is sometimes convenient to set 
\bea
&&\al^{'\vee}_{i,l}=\frac{1-p^{*l}}{1-p^l}q^{lc} \al^\vee_{i,l}. 
\ena
Then we have
\bea
&& [\al^{'\vee}_{i,l},\al^{'\vee}_{j,m}]=\delta_{l+m,0}\frac{[a_{ij} l]_i}{l}\frac{q^{lc}-q^{-lc}}{q_j-q_j^{-1}}\frac{1-p^{*l}} {1-p^l}
\kappa^{-lm_{ij}}q^{lc}.
\ena

Setting also 
 \begin{eqnarray}
&&\alpha_{i,m} =[d_i] \alpha_{i,m}^\vee \qquad 
 \alpha'_{i,m} =
[d_i]\alpha^{'\vee}_{i,m},           \lb{al}     
\end{eqnarray}
one obtains the following commutation relations free from $q_i$\cite{Drinfeld}. 
\begin{lem}\lb{lem:sec4}
For $m,\,n \in \mathbb{Z}_{\ne 0}$, the following commutation relations hold.
\begin{align}
& [\alpha_{i,m},\alpha_{j,n}]=\delta_{m+n,0}\frac{[b_{ij}m]}{m} \frac{q^{cm}-q^{-cm}}{q-q^{-1}}
\frac{1-p^m}{1-p^{*m}}\kappa^{-m\, m_{ij}}q^{-cm}, \\
& [\alpha'_{i,m},\alpha'_{j,n}]=\delta_{m+n,0}\frac{[b_{ij}m]}{m} \frac{q^{cm}-q^{-cm}}{q-q^{-1}}
\frac{1-p^{*m}}{1-p^{m}}\kappa^{-m\, m_{ij}}q^{cm},\\
& [\alpha_{i,m},\alpha'_{j,n}]=[\alpha'_{i,m},\alpha_{j,n}]=\delta_{m+n,0}\frac{[b_{ij}m]}{m}\frac{q^{cm}-q^{-cm}}{q-q^{-1}}\kappa^{-m\,m_{ij}}, \\
& [\alpha_{i,m}, x_j^+(z)]=\frac{[b_{ij}m]}{m}\frac{1-p^m}{1-p^{*m}}q^{-cm}\kappa^{-m\,m_{ij}}z^m x_j^+(z), \\
& [\alpha_{i,m}, x_j^-(z)]=-\frac{[b_{ij}m]}{m}\kappa^{-m\,m_{ij}} z^m x_j^-(z),
\\
& [\alpha'_{i,m}, x_j^+(z)]=\frac{[b_{ij}m]}{m}\kappa^{-m\,m_{ij}}z^m x_j^+(z), \\
& [\alpha'_{i,m}, x_j^-(z)]=-\frac{[b_{ij}m]}{m}\frac{1-p^{*m}}{1-p^m}q^{cm}\kappa^{-m\,m_{ij}} z^m x_j^-(z).
\end{align}
\end{lem}

\begin{prop}
Let $P=\sum_{i}c_iP_i, P+h=\sum_{i}c_i(P_i+h_i) \in H$. 
Define a map $\omega\ :\ U_{q,\kappa,p}(\g) \to U_{q^{-1},\kappa,p^*}(\g)$ by
\be
&&P \mapsto P+h,\quad  P+h \mapsto P,\quad x^+_i(z)\mapsto x^-_i(z),\quad x^-_i(z)\mapsto x^+_i(z),\quad
\al^\vee_{i,k} \mapsto -
\al^{'\vee}_{i,k},\quad\\
&& K^\pm_i \mapsto K^\pm_i,\quad  q^{\pm c/2} \mapsto q^{\mp c/2}, \quad q^{d}\mapsto q^{-d}
\en
for $i\in I$, $k\in \Z_{\not=0}$. Then $\omega$ is an automorphism ( isomorphism). 
\end{prop}
Note that in $U_{q^{-1},\kappa,p^*}(\g)$ we have $(p^*)^*=p^*(q^{-1})^{-2c}=p$.  

\medskip
\noindent
{\it Remark.} \ For representations, on which $q^{\pm c/2}$ take complex values e.g. $q^{\pm k/2}$ ( in Sec.\ref{Zalgstr}, $k=1$ and Sec.\ref{Seclevelzero},  $k=0$), we treat $p, p^*$ as generic complex numbers and assume  $|p|<1$ and $|p^*|=|pq^{-2c}|<1$. 
Then  one has 
\be
 &&g_{ij}(z;s)=\frac{(sq^{b_{ij}}z;s)_\infty}{(sq^{-b_{ij}}z;s)_\infty}
\en
for $|sq^{\pm b_{ij}}z|<1$, $s=p, p^*$, where we set
\be
&&(z;s)_\infty=\prod_{n=0}^\infty(1-zs^n).
\en
Hence in the sense of analytic continuation, one can rewrite \eqref{xpxpgen} and \eqref{xmxmgen} as
\bea
&&z\theta_{p^*}(q^{b_{ij}}\kappa^{-m_{ij}}w/z)x^+_i(z)x^+_j(w)
=-w\kappa^{-m_{ij}}\theta_{p^*}(q^{b_{ij}}\kappa^{m_{ij}}z/w)x^+_j(w)x^+_i(z),\lb{xpxp}\\
&&z\theta_{p}(q^{-b_{ij}}\kappa^{-m_{ij}}w/z)x^-_i(z)x^-_j(w)
=-w\kappa^{-m_{ij}}\theta_{p}(q^{-b_{ij}}\kappa^{m_{ij}}z/w)x^-_j(w)x^-_i(z).\lb{xmxm}
\ena
Here $\theta_p(z)$ denotes Jacobi's odd theta function given by 
\be
\theta_{s}(z)=(z;s)_\infty(s/z;s)_\infty
\en
for $s\in \C, |s|<1$. 
Similarly, one can derive the following relations. 
\begin{prop} 
\bea
&&\phi^\pm_i(z)\phi_j^\pm(w)=\frac{\theta_p(q^{b_{ij}}\kappa^{-m_{i,j}}w/z)}{\theta_p(q^{-b_{ij}}\kappa^{-m_{i,j}} w/z)}
\frac{\theta_{p^*}(q^{-b_{ij}}\kappa^{-m_{i,j}} w/z)}{\theta_{p^*}(q^{b_{ij}}\kappa^{-m_{i,j}} w/z)}
\phi_j^\pm(w)\phi^\pm_i(z),\\
&&\phi^+_i(z)\phi_j^-(w)=\frac{\theta_p(q^{b_{ij}}\kappa^{-m_{i,j}}q^k w/z)}{\theta_p(q^{-b_{ij}}\kappa^{-m_{i,j}}q^k w/z)}
\frac{\theta_{p^*}(q^{-b_{ij}}\kappa^{-m_{i,j}}q^{-k} w/z)}{\theta_{p^*}(q^{b_{ij}}\kappa^{-m_{i,j}}q^{-k}w/z)}
\phi_j^-(w)\phi^+_i(z),\\
&&\phi^\pm_i(z)x^+_j(w)=q^{b_{ij}}\frac{\theta_{p^*}(q^{-b_{ij}}\kappa^{-m_{i,j}}q^{\mp {k}/{2}} w/z)}{\theta_{p^*}(q^{b_{ij}}\kappa^{- m_{i,j}}q^{\mp{k}/{2}}w/z)}x^+_j(w)\phi^\pm_i(z),\lb{phixp}\\
&&\phi^\pm_i(z)x^-_j(w)=q^{-b_{ij}}\frac{\theta_{p}(q^{b_{ij}}\kappa^{-m_{i,j}}q^{\pm {k}/{2}} w/z)}{\theta_{p}(q^{-b_{ij}}\kappa^{- m_{i,j}}q^{\pm{k}/{2}}w/z)}x^-_j(w)\phi^\pm_i(z).\lb{phixm}
\ena
\end{prop}

\subsection{Algebra homomorphism from  $U_{q,\kappa,p}(\g_{tor})$
to $(U_{q,\kappa}(\g_{tor})\otimes \FF[[p]])\sharp \C[\cR_Q]$} \lb{app:1.2}

Let $\C[\cR_Q]$ denote the group algebra of  $\cR_Q$, i.e. $e^{Q_\al}, e^{Q_\beta},  e^{Q_\al}e^{Q_\beta}=e^{Q_\al+Q_\beta}, e^0=1\in \C[\cR_Q]$.
 Consider a composition of algebras $(U_{q,\kappa}(\g_{tor})\otimes \FF[[p]])\sharp \C[\cR_Q]$ with the smash product $\sharp$ defined by
\be
&&g(h,P)a\otimes e^{Q_\al}\cdot f(h,P)b\otimes e^{Q_\beta}\\
&&=g(h,P)f(h-\langle \wt(a),h\rangle,P-\langle Q_\al,P\rangle)ab\otimes e^{Q_\al+Q_\beta}\qquad 
\en
for $a, b \in U_{q,\kappa}(\g_{tor}),\ g(h, P), f(h,P)\in \FF[[p]],\ 
e^{Q_\al}, e^{Q_\beta}\in \C[\cR_Q].$

\begin{prop}
Let us define $u_i^\pm(z,p)\ (i\in I)$ by
\be
u_i^+(z,p)&=&\exp\left(-(q_i-q_i^{-1})\sum_{n>0}\frac{a^\vee_{i,-n}}{1-p^{*n}}(q^{\frac{c}{2}}p^* z)^n\right)\ \in U_{q,\kappa}(\g_{tor})[[p^*]][[z,z^{-1}]] , \\
u_i^-(z,p)&=&\exp\left((q_i-q_i^{-1})\sum_{n>0}\frac{a^\vee_{i,n}}{1-p^n}(q^{\frac{c}{2}} p^{-1} z)^{-n}\right) \ \in U_{q,\kappa}(\g_{tor})[[p]][[z,z^{-1}]].
\en
Then the following gives a homomorphism from  $U_{q,\kappa,p}(\g_{tor})$ to 
$(U_{q,\kappa}(\g_{tor})\otimes \FF[[p]])\sharp \C[\cR_Q]$. 
\bea
&& x^+_i(z)\ \mapsto\ u_i^+(z,p)X_i^+(z)e^{-Q_{i}}, 
\\
&& x^-_i(z)\ \mapsto\ X_i^-(z)u_i^-(z,p),
\\
&&\phi^{\pm}_{i}(z)
\ \mapsto\ u_i^+(q^{\pm\frac{c}{2}} z,p)\Phi^\pm_i(z)u_i^-(q^{\mp\frac{c}{2}}z,p)e^{-Q_{i}}. 
\lb{gdress4}
\ena
\end{prop}
The statement follows from the next Lemma and 
\be
&&[e^{-Q_{i}}, u^\pm_j(z,p)]=0=[e^{-Q_{i}},X^\pm_j(z)].
\en
\begin{lem}\lb{lem:app1}
The following relations hold.
\begin{eqnarray}
&&u_i^+(z,p)X^+_j(w)=
g_{ij}(\kappa^{m_{ij}}z/w;p^*)
X^+_j(w)u_i^+(z,p),\\
&&u_i^+(z,p)X^-_j(w)=g_{ij}(\kappa^{m_{ij}}q^cz/w;p^*)^{-1}
X^-_j(w)u_i^+(z,p),\\
&&u_i^-(z,p)X^+_j(w)=g_{ij}(\kappa^{-m_{ij}}q^{-c}w/z;p)^{-1}
X^+_j(w)u_i^-(z,p),\\
&&u_i^-(z,p)X^-_j(w)=g_{ij}(\kappa^{-m_{ij}}w/z;p)
X^-_j(w)u_i^-(z,p)
,\\
&&\Phi^+_i(z)u^+_j(w,p)=g_{ij}(p^*\kappa^{-m_{ij}}q^{3c/2}w/z;p^*)
g_{ij}(p^*\kappa^{-m_{ij}}q^{-c/2}w/z;p^*)^{-1}u^+_j(w,p)\Phi^+_i(z),\nn\\
&&\\
&&\Phi^-_i(z)u^-_j(w,p)=g_{ij}(p\kappa^{-m_{ij}}q^{-3c/2}z/w;p)
g_{ij}(p\kappa^{-m_{ij}}q^{c/2}z/w;p)^{-1}u^-_j(w,p)\Phi^-_i(z),\\
&&[\Phi^\pm_i(z),u^\mp_j(w,p)]=0,\\
&&u_i^-(z,p)u_j^+(w,p)=g_{ij}(p\kappa^{-m_{ij}}q^{-c}w/z;p)
g_{ij}(p^*\kappa^{-m_{ij}}q^{c}w/z;p^*)^{-1}
u_j^+(w,p)u_i^-(z,p).\nonumber\\
&&
\end{eqnarray}
\end{lem}
This homomorphism is invertible. Hence one has the following isomorphism.
\begin{thm}\lb{isom}
For generic $q, \kappa, p$, 
\be
&&U_{q,\kappa,p}(\g_{tor})\ \cong \ 
(U_{q,\kappa}(\g_{tor})\otimes \FF[[p]])\sharp \C[\cR_Q].
\en
\end{thm}

\subsection{Horizontal and vertical subalgebras}
The elliptic quantum toroidal algebra ${U}_{q,\kappa,p}(\g_{tor})$ has two elliptic quantum algebras $U_{q,p}(\gh)$ as subalgebras,  
which we call the horizontal and vertical elliptic quantum algebras. This is an elliptic analogue of Proposition \ref{propGKV}. 

\begin{prop}
The elliptic quantum toroidal algebra ${U}_{q,\kappa,p}(\g_{tor})$ has the elliptic quantum algebra  
${U}_{q,p}(\gh)$ generated by $\al^\vee_{i,l}, k_i^\pm, x^+_{i,n}, x^-_{i,n}, q^{\pm {c}/{2}}$ $(i\in \ I\backslash \{0\}, l\in \Z\backslash \{0\}, n\in \Z  
)$ as a subalgebra. We call it the horizontal subalgebra.
The ${U}_{q,\kappa,p}(\g_{tor})$ has another subalgebra  
${U}_{q,p'}(\gh)$ generated by $x^+_{i,0}, x^-_{i,0},  K^{\pm}_i$ $(i\in I)$, which is again the elliptic quantum algebra with $p'$ being not necessalily equal to $p$.  
When we take $p'=p$, we call it the vertical subalgebra. 
\end{prop}
\noindent
{\it Proof.} The first half is obvious by comparing the relations in Definition \ref{defEQTA} and those of $U_{q,p}(\gh)$ for example see \cite{FKO}.  
 The second half  is due to the fact that 
$x^+_{i,0}, x^-_{i,0},  K^\pm_i, q^{\pm c/2}$ $(i\in I)$ 
generate the dynamical quantum affine algebra $(\FF\otimes_\C U_{q}(\gh))\sharp\C[\cR_Q]$ in the Drinfeld-Jimbo formulation\cite{Drinfeld85,Jimbo85}.  
In fact from Definition \ref{defEQTA} one can find the following relations. 
\bea
&&
g(h,{P})x^+_{j,0}=x^+_{j,0}\ g(h+\langle \al_{j},h\rangle, P-\langle Q_{j},P\rangle),\lb{ge0}\\
&&
g(h,{P})x^-_{j,0}=x^-_{j,0}\ g(h-\langle \al_{j},h\rangle, P),\lb{gf0}\\
&&g(h,{P})K^{\pm}_j=K^{\pm}_j\ g(h,P-\langle Q_{j},P\rangle),
\lb{gKpm0}\\
&&K^\pm_ix^+_{j,0}(K^\pm_i)^{-1}=q^{\pm b_{ij}}x^+_{j,0}, \lb{Kixp0}\\
&&K^\pm_ix^-_{j,0}(K^\pm_i)^{-1}=q^{\mp b_{ij}}x^-_{j,0},\\
&&[x^+_{i,0},x^-_{j,0}]=\frac{\delta_{i,j}}{q_i-q_i^{-1}}(K^+_i-K^-_i),\\
&&\sum_{r=0}^a(-1)^r\left[\mmatrix{a\cr r\cr}\right]_i
(x^\pm_{i,0})^{r}x^\pm_{j,0} (x^\pm_{i,0})^{a-r}=0\qquad \mbox{for}\ i\not=j.\lb{xpm0Serre}
\ena
Hence by the isomorphism\cite{Beck} between the  Drinfeld-Jimbo\cite{Drinfeld85,Jimbo85}  and the Drinfeld\cite{Drinfeld} formulations of $U_q(\gh)$, one can construct the Drinfeld generators $a^\vee_{i,l}, X^+_{i,n}, X^-_{i,n} \ (i\in I\backslash \{0\}, \ l \in \Z\backslash \{0\}, \ 
n \in \Z)$ for the latter in terms of $x^+_{i,0}, x^-_{i,0},  K^{\pm}_i$. Then applying the isomorphism in Theorem \ref{isom} with using $p'$ as $p$, one obtains the elliptic quantum algebra ${U}_{q,p'}(\gh) $.    
\qed

\section{$H$-Hopf Algebroid Structure}
We introduce a $H$-Hopf algebroid structure\cite{EV98,EV99,KR,Konno09} to $U_{q,\kappa,p}(\g_{tor})$ as its co-algebra 
structure. 

\subsection{$H$-algebra}
Let $\cA$ be an associative algebra, $\cH$ be a finite dimensional subspace of 
$\cA$, and $\cH$ be its dual. Denote by $\cM_{\cH^*}$  the 
field of meromorphic functions on $\cH^*$. 

\begin{dfn}[$\cH$-algebra] 
A $\cH$-algebra is an associative algebra $\cA$ with 1, which is bigraded over 
$\cH^*$, $\ds{\cA=\bigoplus_{\alpha,\beta\in \cH^*} \cA_{\al\beta}}$, and equipped with two 
algebra embeddings $\mu_l, \mu_r : \cM_{\cH^*}\to \cA_{00}$ (the left and right moment maps), such that 
\be
\mu_l(\hf)a=a \mu_l(T_\al \hf), \quad \mu_r(\hf)a=a \mu_r(T_\beta \hf), \qquad 
a\in \cA_{\al\beta},\ \hf\in \cM_{\cH^*},
\en
where $T_\al$ denotes a difference operator $(T_\al \hf)(\la)=\hf(\la+\al)$ of $\cM_{\cH^*}$.
\end{dfn}

Let $\cA$ and $\cB$ be two $\cH$-algebras. 
\begin{dfn}
A $\cH$-algebra homomorphism from $\cA$ to $\cB$ is an algebra homomorphism $\pi:\cA\to \cB$ 
preserving the bigrading and the moment maps, i.e. $\pi(\cA_{\al\beta})\subseteq \cB_{\al\beta}$ for all $\al,\beta \in \cH^*$ and $\pi(\mu^\cA_l(\hf))=\mu^\cB_l(\hf), \pi(\mu^\cA_r(\hf))=\mu^\cB_r(\hf)$. 
\end{dfn}

The tensor product $\cA {\widetilde{\otimes}}\cB$ is the $\cH^*$-bigraded vector space with 
\be
 (\cA {\widetilde{\otimes}}\cB)_{\al\beta}=\bigoplus_{\gamma\in\H^*} (\cA_{\al\gamma}\otimes_{\cM_{\cH^*}}\cB_{\gamma\beta}),
\en
where $\otimes_{\cM_{\cH^*}}$ denotes the usual tensor product 
modulo the following 
relation.
\bea
\mu_r^\cA(\hf) a\otimes b=a\otimes\mu_l^\cB(\hf) b, \qquad a\in \cA, 
b\in \cB, \hf\in \cM_{\cH^*}.\lb{AtotB}
\ena
The tensor product $\cA {\widetilde{\otimes}}\cB$ is again a $\cH$-algebra with the multiplication $(a\otimes b)(c\otimes d)=ac\otimes bd$ and the moment maps 
\be
\mu_l^{\cA {\widetilde{\otimes}}\cB} =\mu_l^\cA\otimes 1,\qquad \mu_r^{\cA {\widetilde{\otimes}}\cB} =1\otimes \mu_r^\cB.
\en

Let $\cD$ be the algebra of difference operators on $\cM_{\cH^*}$
\be
\cD&=&\{\ \sum_i \hf_i T_{\beta_i}\ |\ \hf_i\in \cM_{\cH^*},\ \beta_i\in \cH^*\ \}.
\en
Equipped  with the bigrading 
 $\cD_{\al\al}=\{\  \hf T_{-\al}\ |\ \hf\in \cM_{\cH^*},\ \al\in \cH^*\ \}$, 
 $\cD_{\al\beta}=0\ (\al\not=\beta)$ 
 and the moment maps $\mu^{\cD}_l, \mu^{\cD}_r : \cM_{\cH^*}\to \cD_{00}$ 
 defined by 
$\mu^{\cD}_l(\hf)=\mu^{\cD}_r(\hf)=\hf T_0$, $\cD$ is a $\cH$-algebra.
 For any $\cH$-algebra $\cA$, we have the canonical 
isomorphism as an $\cH$-algebra 
\bea
&&\cA\cong \cA\tot \cD\cong  \cD\tot \cA \lb{Diso} 
\ena
by $a\cong a\tot T_{-\beta}\cong T_{-\al}\tot a$ for all $a\in \cA_{\al\beta}$.

\begin{prop}
Let $P+h=\sum_{i\in I}c_i(P_i+h_i)\in H,\ c_i\in \C$. 
$\cU=U_{q,\kappa,p}(\gt)$  is {an $H$-algebra} by 
\be
&&\cU=\bigoplus_{\al,\beta\in \h^*}\cU_{\al,\beta}\\
&&(\cU)_{\al\beta}=\left\{x\in U \left|\ q^{P+h}x q^{-(P+h)}=q^{\langle\al,P+h\rangle}x,\quad q^{P}x q^{-P}=q^{\langle Q_\beta,P\rangle}x\quad \forall P+h, P\in H\right.\right\}
\en
and $\mu_l, \mu_r : \FF \to \cU_{0,0}$ defined by 
\be
&&\mu_l(\hf)=f(h,P+h,p)\in \FF[[p]],\qquad \mu_r(\hf)=f(h,P,p^*)\in \FF[[p^*]]
\en
for $\hf=f(h,P,p^*)\in \FF[[p^*]]$.
\end{prop}
Corresponding to \eqref{AtotB}, we have 
\bea
&&f(h,P,p^*)a\tot b=a\tot f(h,P+h,p)b\qquad a,b\in \cU,\lb{fsabafb}.
\ena
Note that $p=p^*q^{2c}$.
We also 
regard $T_\al=e^{-Q_\al}\in \C[\cR_Q]$ as a shift operator 
\be
&&(T_\al \widehat{f})=f(h,P+\langle Q_\al,P\rangle,p)
\en
Then $\cD=\{\widehat{f} e^{-Q_\al} \ |\  \widehat{f}\in \FF,  e^{-Q_\al}\in \C[\cR_Q]\}$ becomes the $H$-algebra  having the property \eqref{Diso} 
for $\cA=\cU$.

\subsection{Hopf algebroid}
\begin{dfn}\lb{Deltavep}
A $\cH$-bialgebroid is a $\cH$-algebra $\cA$ equipped with two $\cH$-algebra homomorphisms 
$\Delta:\cA\to \cA{{\tot}}\cA$ (the comultiplication) and $\vep : \cA\to \cD$ (the counit) such that 
\be
&&(\Delta \tot \id)\circ \Delta=(\id \tot \Delta)\circ \Delta,\\
&&(\vep \tot \id)\circ\Delta =\id =(\id \tot \vep)\circ \Delta,
\en
under the identification \eqref{Diso}.
\end{dfn}
 
\begin{dfn}
\lb{defS}
A $\cH$-Hopf algebroid is a $\cH$-bialgebroid $\cA$ equipped with a $\C$-linear map $S : \cA\to \cA$ (the antipode), such that 
\be
&&S(\mu_r(\hf)a)=S(a)\mu_l(\hf),\quad S(a\mu_l(\hf))=\mu_r(\hf)S(a),\quad \forall a\in \cA, \hf\in \cM_{\cH^*},\\
&&m\circ (\id \tot S)\circ\Delta(a)=\mu_l(\vep(a)1),\quad \forall a\in \cA,\\
&&m\circ (S\tot\id  )\circ\Delta(a)=\mu_r(T_{\al}(\vep(a)1)),\quad \forall a\in \cA_{\al\beta},
\en
where $m : \cA{{\tot}} \cA \to \cA$ denotes the multiplication and $\vep(a)1$ is the result of applying the difference operator $\vep(a)$ to the constant function $1\in \cM_{\cH^*}$.
\end{dfn}

Note that the $\cH$-algebra $\cD$ is a $\cH$-Hopf algebroid with 
$\Delta_\cD : \cD\to \cD\tot \cD,\ \vep_\cD: \cD \to \cD,\ 
S_\cD : \cD \to \cD$ defined by 
\be
&&\Delta_\cD(\hf T_{-\al})=\hf T_{-\al} \tot T_{-\al},\\
&&\vep_\cD=\id,
\qquad  S_\cD(\hf T_{-\al})=T_{\al}\hf=(T_{\al}\hf)T_{\al}.
\en

Now let us consider our $H$-algebra $\cU$. 
We define two $H$-algebra homomorphisms, the co-unit $\vep : \cU\to \cD$ and the co-multiplication $\Delta : \cU\to \cU\widetilde{\otimes}\cU$ by
\bea
&&\vep(x^\pm_i(z))=0,\quad \vep(\phi^\pm_i(z))=T_{\al_i},\\
&&\vep(q^{\pm c/2})=1, \quad \vep(e^{Q_\al})=e^{Q_\al}\quad (Q_\al\in \cR),\\
&&\vep(\mu_l({\hf}))= \vep(\mu_r(\hf))=\widehat{f}T_0, \lb{counitf}
\ena
\bea
&&\Delta(q^{\pm c/2})=q^{\pm c/2}\tot q^{\pm c/2},\\
&&\Delta(\phi^\pm_i(z))=\phi^\pm_i(q^{\mp c^{(2)}/2}z)\tot \phi^\pm_i(q^{\pm c^{(1)}/2}z),\lb{co1}\\
&&\Delta(x^+_i(z))=1\tot x^+_i(z)+x^+_i(q^{c^{(2)}}z)\tot\phi^-_i(q^{c/2}z),\\
&&\Delta(x^-_i(z))= x^-_i(z)\tot 1+\phi^+_i(q^{c/2}z)\tot x^-_i(q^{c^{(1)}}z),\lb{co3}\\
&&\Delta(e^{Q_\al})=e^{Q_\al}\tot e^{Q_\al},\\
&&\Delta(\mu_l(\hf))=\mu_l(\hf)\tot 1,\qquad \Delta(\mu_r(\hf))=1\tot \mu_r(\hf).
\ena
Here $q^{\pm c^{(1)}/2}=q^{\pm c/2}\tot 1, \ q^{\pm c^{(2)}/2}=1\tot q^{\pm c/2} $.
One can easily check that these satisfy the relations in Definition \ref{Deltavep}. 
This $\Delta$ is called the Drinfeld comultiplication. 

Define further an algebra antihomomorphism (the antipode) $S : \cU\to \cU$ by
\bea
&&S(x^+_{i}(z))=-x^+_i(q^{-c}z)\phi^-_i(q^{-c/2}z)^{-1},\lb{Sonxp}\\
&&S(x^-_{i}(z))=-\phi^+_i(q^{-c/2}z)^{-1}x^-_i(q^{-c}z),\lb{Sonxm}\\
&&S(q^{\pm c/2})=q^{\mp c/2},\\
&&S(e^{Q_\al})=e^{-Q_{\al}},\quad S(\mu_r(\hf))=\mu_l(\hf),\quad S(\mu_l(\hf))=\mu_r(\hf).
\ena

\begin{thm}
The $H$-algebra $\cU$ equipped with $(\Delta,\vep,S)$ is an $H$-Hopf algebroid. 
\end{thm}

\begin{dfn}
We call the $H$-Hopf algebroid $(\cU,H,{\cM}_{H^*},\mu_l,\mu_r,\Delta,\vep,S)$ the  toroidal  elliptic quantum  group $U_{q,\kappa,p}(\gt)$. 
\end{dfn}

\noindent
{\it Remark.}\ 
In a representation where $q^{\pm c/2}=1$ such as the level-$(0,l)$ representation  in Sec.\ref{Seclevelzero}, the following gives a comultiplication, too.
 \bea
&&\Delta^{op}(\phi^\pm_i(z))=\phi^\pm_i(z)\tot \phi^\pm_i(z),\lb{opco1}\\
&&\Delta^{op}(x^+_i(z))=x^+_i(z) \tot 1+\phi^-_i(z)\tot x^+_i(z),\\
&&\Delta^{op}(x^-_i(z))= 1\tot x^-_i(z)+ x^-_i(z)\tot \phi^+_i(z).\lb{opco3}
\ena

\section{Dynamical Representations}
Let us consider a vector space $\hV$ over $\FF$, which is  
${H}$-diagonalizable, i.e.  
\be
&&\hV=\bigoplus_{\la,\mu\in {\h}^*}\hV_{\la,\mu},\ \hV_{\la,\mu}=\{ v\in \hV\ |\ q^{P+h}\cdot v=q^{\langle\la,P+h\rangle} v,\ q^{P}\cdot v=q^{\langle Q_\mu,P\rangle} v,\ \forall 
P+h, P\in 
{H}\}.
\en
Let us define the $H$-algebra $\cD_{H,\hV}$ of the $\C$-linear operators on $\hV$ by
\be
&&\cD_{H,\hV}=\bigoplus_{\al,\beta\in {\h}^*}(\cD_{H,\hV})_{\al\beta},\\
&&\hspace*{-10mm}(\cD_{H,\hV})_{\al\beta}=
\left\{\ X\in \End_{\C}\hV\ \left|\ 
\mmatrix{ q^{P+h}X q^{-(P+h)}=q^{\langle\alpha,P+h\rangle} X ,
\quad q^P Xq^{-P}=q^{\langle Q_\beta,P\rangle}X{,}\cr
 X\cdot\hV_{\la,\mu}\subseteq 
 \hV_{\la+\al,\mu+\beta}\cr}  
 \right.\right\},\\
&&\mu_l^{\cD_{H,\hV}}(\widehat{f})v=f(\langle\la,h\rangle, \langle\la,P+h\rangle,p)v,\quad 
\mu_r^{\cD_{H,\hV}}(\widehat{f})v=f(\langle\la,h\rangle, \langle Q_\mu,P\rangle,p^*)v,\quad \widehat{f}\in {\cM}_{H^*},
\ v\in \hV_{\la,\mu}.
\en
\begin{dfn}
We define a dynamical representation of $U_{q,\kappa,p}(\gt)$ on $\hV$ to be  
 a $H$-algebra homomorphism ${\pi}: U_{q,\kappa,p}(\gt) 
 \to \cD_{H,\hV}$. 
\end{dfn}
\begin{dfn}
For $k, l\in \C$, we say that a $U_{q,\kappa,p}(\gt)$-module has  level $(k,l)$ if the two central elements $q^{c/2}$ and $\prod_{i\in I}(K_i^+)^{a^\vee_i}$ act 
as $q^{k/2}$ and $q^{-l}$ on it, respectively.  
\end{dfn}

\begin{dfn}
For $\omega\in \C$, we set 
\be
&&\hV_\omega=\{v\in \hV\ |\ q^{-d}\cdot v=q^\omega v\ \}
\en
and we call $\hV_\omega$ the space of elements homogeneous of degree $\omega$. 
We also say that $X\in \cD_{H,\hV}$ is homogeneous of degree $\omega\in \C$ if 
\be
&&[q^{-d}, X]=q^\omega X
\en
and denote by $(\cD_{H,\hV})_\omega$ the space of all {endomorphisms} homogeneous of degree $\omega$. 
\end{dfn}

\begin{dfn}
Let $U({\frak H})$, $U({\frak N}_+), U({\frak N}_-)$ be the subalgebras of 
$U_{q,\kappa,p}(\gt)$ 
generated by 
$q^{\pm c/2}, q^d, 
K^\pm_{i}\ (i\in I)$, by $\al^\vee_{i,n}\ (i\in I, n\in \Z_{>0})$,  $x^+_{i, n}\ (i\in I, n\in \Z_{\geq 0})$  
$x^-_{i, n}\ (i\in I, n\in \Z_{>0})$ and by $\al^\vee_{i,-n}\ (i\in I, n\in \Z_{>0}),\ x^+_{i, -n}\ (i\in I, n\in \Z_{> 0}),\ x^-_{i, -n}\ (i\in I, n\in \Z_{\geq 0})$, respectively.   
\end{dfn}

\begin{dfn}
For $k, l\in\C$, $\la, \mu\in \h^*$, 
a (dynamical) $U_{q,\kappa,p}(\gt)$-module $\hV(\la,\mu)$ is called the 
level-$(k,l)$ highest weight module with the highest weight $(\la,\mu)$, if there exists a vector 
$v\in \hV(\la,\mu)$ such that
\be
&&\hV(\la,\mu)=U_{q,\kappa,p}(\gt)\cdot v,\qquad U(\gN_+)\cdot v=0,\\
&&q^{c/2}\cdot v=q^{k/2}v,\qquad  \prod_{i\in I}(K_i^+)^{a^\vee_i}\cdot v
=q^{-l}v \quad \\
&& f(h,{P})\cdot v =f(\langle \la-\mu, h\rangle, {\langle Q_\mu,P\rangle})v,
\en
\end{dfn}
Note that 
\be
&&q^P\cdot v=q^{\langle Q_\mu,P\rangle}v,\quad q^{P+h}\cdot v=q^{\langle \la-\mu+Q_\mu,P+h\rangle}v=q^{\langle \la,P+h\rangle}v.
\en
Hence $\la, \mu$ are weights measured by $P+h$ and $P$, respectively.  

We define the category $\gC_{k,l}$ as follows. 
\begin{dfn}
For $k, l\in \C$, $\gC_{k,l}$ is the full subcategory of the category of $U_{q,\kappa,p}(\gt)$-modules 
consisting of those modules $\hV$ such that 
\begin{itemize}
\item[(i)] $\hV$ has level $(k,l)$
\item[(ii)] $\ds{\hV=\bigsqcup_{\omega\in \C} \hV_\omega}$
\item[(iii)] For every $\omega \in \C$, there exists $n_0\in \N$ such that for all 
$n>n_0$, $\hV_{\omega+n}=0$.   
\end{itemize}
\end{dfn}

Since $\hpi U(\gN_+)\subset \bigsqcup_{n\in \Z_{\geq 0}}(\cD_{H,\hV})_n$, any level-$(k,l)$ highest weight $U_{q,\kappa,p}(\gt)$-modules belong to $\gC_{k,l}$.

\section{$Z$-algebra Structure of $U_{q,\kappa,p}(\gt)$}\lb{Zalgstr}
In this section we discuss the $Z$-algebra structure for the elliptic quantum toroidal algebras $U_{q,\kappa,p}(\gt)$. We show that the $Z$-algebra of $U_{q,\kappa,p}(\gt)$ is essentially  the same as the one of the quantum toroidal algebra $U_{q,\kappa}(\gt)$ except for that  the former is a dynamical version of the latter. 
 In both cases the  $Z$-algebras govern the irreducibility of the level-$(k(\not=0), l)$ $U_{q,\kappa,p}$ or $U_{q,\kappa}$-modules.  As an example,  we construct the level $(1,l)$ irreducible representation of $U_{q,\kappa,p}(\g_{tor})$ for the simply laced $\g_{tor}$. These results are toroidal analogues of those for the elliptic quantum algebras $U_{q,p}({\gh})$ and  the quantum affine algebras $U_{q}({\gh})$ discussed in \cite{FKO}.  The $Z$-algebra was first introduced  in the representation theory of affine Lie algebras 
 by Lepowsky and Wilson\cite{LW}, and its $q$-analogue associated with $U_{q}({\gh})$ was investigated 
 in \cite{BoVi,Jing96,Jing00}.

\subsection{Heisenberg subalgebra $U(\gH)$}\lb{Heisenberg}
Let 
$U(\gH_+)$ (resp. $U({\gH}_-)$) be the commutative
 subalgebras of  $U({{\gH}})$ generated by $\{q^{\pm c/2},\ \prod_{i\in I}(K_i^\pm)^{a^\vee_i},\\ \alpha_{i,n}(i\in I, n\in\Z_{> 0})\}$ 
(resp. $\{\alpha_{i,-n}(i\in I, n\in\Z_{> 0})\}$). 
We have
\be
U(\gH)=U(\gH_-)U(\gH_+).
\en

For $k,l\in \C$, let  $\C 1_{k,l}$ be the one-dimensional 
 $U(\gH_+)$-module generated by the vacuum vector $1_{k,l}$ defined by  
\be
&&q^{c/2}\cdot 1_{k,l}=q^{k/2}1_{k,l} 
\qquad \prod_{i\in I}(K_i^+)^{a^\vee_i}\cdot 1_{k,l}=q^{-l}1_{k,l},  
\qquad \alpha_{i,n}\cdot 1_{k,l}=0 \qquad(n> 0). 
\en
Then we have the induced $U({\gH})$-module
\be
\F^{\al}_{k,l}=U(\gH)\otimes_{U({\gH}_+)}\C1_{k,l}.
\en
The space  $\F^{\al}_{k,l}$ can be identified with a polynomial ring $\C[\al_{i,-m}\ (i\in I,m>0)]$ 
by  
\be
&&q^{c/2}\cdot u= q^{k/2}u,\qquad \prod_{i\in I}(K_i^+)^{a^\vee_i}\cdot u=q^{-l}u, \qquad \\
&&\alpha_{i,-n}\cdot u= \alpha_{i,-n}u,\qquad \\
&&\alpha_{i,n}\cdot u=\sum_{j} \frac{[b_{ij}n][kn]}{n}\frac{1-p^n}{1-p^{*n}}q^{-kn}\frac{\partial}{\partial\alpha_{j,-n}}u\qquad(n>0)
\en
for $u\in \C[\al_{i,-m}\ (i\in I,m>0)]$. We call  $\F^{\al}_{k,l}$ the level-$(k,l)$ Fock module.  
Note that $p^*=pq^{-2k}$ on $\F^{\al}_{k,l}$.

\subsection{Dynamical quantum $Z$-algebra $\cZ_\cV$}

Let $k\in \C^{\times},\ l\in \C$ and consider $(\cV,\pi)\in \gC_{k,l}$. 
We  assume $|p|, |p^*|, |q^{2k}|<1$. Define 
\begin{align}
& E^{\pm}(\alpha_i,z)=\exp\left\{\pm \sum_{m=1}^{\infty}\frac{q-q^{-1}}{q^{km}-q^{-km}} \alpha_{i, \pm m} z^{\mp m} \right\},\\
& E^{\pm}(\alpha'_i,z)=\exp\left\{\mp \sum_{m=1}^{\infty}\frac{q-q^{-1}}{q^{km}-q^{-km}} \alpha'_{i, \pm m} z^{\mp m} \right\}.
\end{align}
One finds these satisfy the relations in Lemma \ref{lem:sec4-2}.
\begin{dfn}
 We define  $\mathcal{Z}_i^{\pm}(z,\cV)=\sum_{m\in \Z} \cZ_{i,m}^\pm(\cV)z^{-m}\in \cD_{H,\cV}[[z,z^{-1}]],\ i \in I$ by
\begin{align}
& \mathcal{Z}_i^+(z,\cV) = E^-(\alpha_i,z)x_i^+(z)E^+(\alpha_i,z), \\
& \mathcal{Z}_i^-(z,\cV) = E^-(\alpha'_i,z)x_i^-(z)E^+(\alpha'_i,z).
\end{align}
\end{dfn}
Due to the truncation property of the grading of $\cV$ w.r.t $q^{-d}$, $\mathcal{Z}_i^{\pm}(z,\cV)$ are well defined operators on $\cV$.

From the defining relations of $U_{q,\kappa,p}(\gt)$ and Lemma \ref{lem:sec4-2}, we obtain the following relations.

\begin{thm}\lb{zz}
The operators $\cZ^+_i(z;\hV)$ satisfy the following relations. For $g(h,P) 
\in \FF$, 
\bea
&& [\alpha_{i,m}, \mathcal{Z}_j^{\pm}(z,\cV)]=0, \quad (m \in \mathbb{Z}_{\ne 0}) \lb{zalg1}\\
&&[\ \prod_{i\in I}(K_i^{\pm} )^{a_i^\vee}, \cZ^+_j(z;\hV)\ ]=0,\lb{KZZK}\\
&&
g(h,{P})\cZ^+_i(z;\hV)=\cZ^+_i(z;\hV)g(h+\langle \al_{i},h\rangle, P-\langle Q_{i},P\rangle),\quad
\\
&&g(h,{P})\cZ^-_i(z;\hV)=\cZ^-_i(z;\hV)g(h-\langle \al_{i},h\rangle, P),\lb{gPZ}\\
&&g(h,{P})K^{\pm}_i=K^{\pm}_i g(h,P\mp\langle Q_{i},P\rangle),\lb{gPK}\\
&&K^\pm_i\cZ_j^+(z,\cV)(K^\pm_i)^{-1}=q^{\pm\langle\al_j,h_i\rangle}\cZ_j^+(z,\cV),\quad 
K^\pm_i\cZ_j^-(z,\cV)(K^\pm_i)^{-1}=q^{\mp\langle\al_j,h_i\rangle}\cZ_j^-(z,\cV),\\
&& {z}\frac{(q^{-b_{ij}}\kappa^{-m_{ij}}\frac{w}{z};q^{2k})_\infty}{(q^{2k}q^{b_{ij}}\kappa^{-m_{ij}} \frac{w}{z};q^{2k})_\infty} \mathcal{Z}_i^{\pm}(z,\cV) \mathcal{Z}_j^{\pm}(w,\cV)=
-{w}\kappa^{-m_{ij}} \frac{(q^{-b_{ij}}\kappa^{m_{ij}}\frac{z}{w};q^{2k})_\infty}{(q^{2k} q^{b_{ij}}\kappa^{m_{ij}} \frac{z}{w};q^{2k})_\infty} 
\mathcal{Z}_j^{\pm}(w,\cV) \mathcal{Z}_i^{\pm}(z,\cV), \nn\\
&&\lb{zalg2}\\
&& \frac{(q^{b_{ij}} q^k \kappa^{-m_{ij}} \frac{w}{z};q^{2k})_{\infty}}{(q^{-b_{ij}}q^k \kappa^{m_{ij}} \frac{w}{z};q^{2k})_{\infty}} \mathcal{Z}_i^+(z,\cV) \mathcal{Z}_j^-(w,\cV) - \frac{(q^{b_{ij}}q^k \kappa^{m_{ij}} \frac{z}{w};q^{2k})_{\infty}}{(q^{-b_{ij}}q^k \kappa^{m_{ij}} \frac{z}{w};q^{2k})_{\infty}} \mathcal{Z}_j^-(w,\cV)\mathcal{Z}_i^+(z,\cV) \nn \\
&& \qquad =\frac{\delta_{i,j}}{q_i-q_i^{-1}} 
\left\{ \delta(q^k \frac{w}{z}) K_i^+ -\delta(q^{-k} \frac{w}{z}) K_i^- \right\}, \lb{zalg3}
\ena
\bea
&&  \sum_{\sigma \in \gS_a} \prod_{1 \leq k < m \leq a} 
\frac{(q^{-2} \frac{z_{\sigma(m)}}{z_{\sigma(k)}};q^{2k})_{\infty}}{(q^{2} \frac{z_{\sigma(m)}}{z_{\sigma(k)}};q^{2k})_{\infty}} \nn\\
&& \qquad \times \sum_{r=0}^{a} (-1)^r \left[ \mmatrix{a \cr r \cr}\right]_i \prod_{1 \leq s \leq r} 
\frac{(q^{-b_{ij}}\kappa^{-m_{ij}}\frac{w}{z_{\sigma(s)}};q^{2k})_{\infty}}{(q^{b_{ij}}\kappa^{-m_{ij}}\frac{w}{z_{\sigma(s)}};q^{2k})_{\infty}}
\prod_{r+1 \leq s \leq a} 
\frac{(q^{-b_{ij}}\kappa^{m_{ij}}\frac{z_{\sigma(s)}}{w};q^{2k})_{\infty}}{(q^{b_{ij}}\kappa^{m_{ij}}\frac{z_{\sigma(s)}}{w};q^{2k})_{\infty}} \nn \\
&& \quad \quad \times \mathcal{Z}_i^+(z_{\sigma (1)},\cV) \cdots \mathcal{Z}_i^+(z_{\sigma(s)},\cV)
\mathcal{Z}_j^+(w,\cV) \mathcal{Z}_i^+(z_{\sigma(s+1)},\cV) \cdots \mathcal{Z}_i^+(z_{\sigma(a)},\cV)=0,  \lb{zalg4}
\\
&&  \sum_{\sigma \in \gS_a} \prod_{1 \leq  k < m \leq a} 
\frac{(q^{-2} q^{2k} \frac{z_{\sigma(m)}}{z_{\sigma(k)}};q^{2k})_{\infty}}{(q^{2} q^{2k} \frac{z_{\sigma(m)}}{z_{\sigma(k)}};q^{2k})_{\infty}} \nn\\
&& \quad \quad \times \sum_{s=0}^{a} (-1)^s \left[ \mmatrix{a \cr r \cr}\right]_i \prod_{1 \leq s \leq r} 
\frac{(q^{-b_{ij}} q^{2k} \kappa^{-m_{ij}} \frac{w}{z_{\sigma(s)}};q^{2k})_{\infty}}{( q^{b_{ij}}q^{2k} \kappa^{-m_{ij}} \frac{w}{z_{\sigma(s)}};q^{2k})_{\infty}}
\prod_{r+1 \leq s \leq a} 
\frac{(q^{-b_{ij}} q^{2k} \kappa^{m_{ij}} \frac{z_{\sigma(s)}}{w};q^{2k})_{\infty}}{(q^{b_{ij}}q^{2k} \kappa^{m_{ij}} \frac{z_{\sigma(s)}}{w};q^{2k})_{\infty}} \nn \\
&& \quad \quad \times \mathcal{Z}_i^-(z_{\sigma (1)},\cV) \cdots \mathcal{Z}_i^-(z_{\sigma(s)},\cV)
\mathcal{Z}_j^-(w,\cV) \mathcal{Z}_i^-(z_{\sigma(s+1)},\cV) \cdots \mathcal{Z}_i^-(z_{\sigma(a)},\cV)=0. \lb{zalg5}
\ena
\end{thm}

\begin{dfn}\label{z-algebra}
For $k\in \C^{\times},\ l\in \C$ and $(\hV,\pi)\in \gC_{k,l}$, we call 
the $H$-subalgebra of $\cD_{H,\hV}$
 generated by $\cZ_{i,m}^\pm(\hV)$, $K^{\pm}_i\ (i\in I, m\in \Z)$ and $q^d$ 
the dynamical quantum $Z$-algebra ${\cZ}_{\hV}$ 
 associated with  $(\hV,\pi)$. 
\end{dfn}

\subsection{The universal $Z$-algebra $\cZ_{k,l}$}
For $k\in \C^{\times}, \ l\in \C$, 
we define the universal dynamical quantum $Z$-algebra $\cZ_{k,l}$ as follows.

\begin{dfn}\label{univz-algebra}
Let  $\cZ^\pm_{i,m}
\ (i\in I, m\in \Z)$ be abstract symbols.  
We set $\cZ^\pm_i(z)=\sum_{m\in \Z}\cZ^\pm_{i,m}z^{-m}$. 
We define  the universal dynamical quantum $Z$-algebra $\cZ_{k,l}$ to be a topological algebra over $\FF[[q^{2k}]]$ generated 
by $\cZ^\pm_{i,m}, K^{\pm }_i \ (i\in I, m\in \Z),  q^d$ 
subject to the relations obtained by replacing 
$\cZ^\pm_i(z;\hV)$ by $\cZ^\pm_i(z)$
 in \eqref{KZZK}-\eqref{zalg5}. 
\end{dfn}
In this Definition, we treat the relations as formal Laurent series in $z, w$ and  $z_j$'s in a similar way to those of the defining relations of $U_{q,\kappa,p}(\gt)$. These relations are well-defined in the $q^{2k}$-adic topology.

\begin{prop}
${\cZ}$ is a $H$-algebra with the same $\mu_l,\mu_r$ as in $U_{q,\kappa,p}(\gt)$. 
\end{prop}

For $(\hV,\hpi)\in \gC_{k,l}, k\in \C^{\times}, \ l\in \C$, we extend $\hpi$ to the map $\hpi: {\cZ}_{k,l}\to \cD_{H,\hV}$ by 
$\hpi(\cZ^\pm_{i,m})=\cZ^\pm_{i,m}(\hV) $. 
Then $\hV$ is a ${\cZ}_{k,l}$-module by $\hpi$.

\begin{dfn}
For $k\in \C^{\times},\ l\in \C$, we denote by $\gD_{k,l}$ the full subcategory of the category of 
${\cZ}_{k,l}$-modules consisting of those modules $(\cW,\sigma)$ such that
\begin{itemize}
\item[(i)] $\cW$ has level $(k,l)$, i.e. $q^{\pm c/2}$ and $\prod_{i\in I}(K_i^{+} )^{a_i^\vee}$ take their values $q^{\pm k/2}$ and $q^{-l}$, respectively, on $\cW$ .
\item[(ii)] $\cW=\bigsqcup_{\omega\in \C} \cW_\omega$, where  $\cW_\omega=\{w\in \cW\ |\ q^{-\sigma(d)} w=q^\omega w\ \}$
\item[(iii)] For every $\omega \in \C$, there exists $n_0\in \N$ such that for all 
$n>n_0$, $\cW_{\omega+n}=0$.   
\end{itemize}
\end{dfn}

For $(\hV,\hpi)\in \gC_{k,l}$, let us define the vacuum space $\Omega_\cV$ as an invariant subspace of the $\cZ_\cV$-module\cite{LW,FKO}. 
\begin{equation}
\Omega_\cV=\{v \in \cV\ |\ \pi(\alpha_{i,m})v=0\quad \forall i\in I, m\in \Z_{>0} \}.
\end{equation}
For a morphism $\varphi\ :\ \cV\to \cV'$ in $\gC_{k,l}$ we have 
\be
&&\varphi(\Omega_{\cV})\subset \Omega_{\cV'}.
\en

\begin{prop}
For $(\hV,\hpi)\in \gC_{k,l}$, there is a unique representation $\sigma$ of ${\cZ}_{k,l}$ on $\Omega_{\hV}$ 
such that $(\Omega_{\hV},\sigma)\in \gD_{k,l}$, 
\be
&&\sigma(K^\pm_i)=\hpi(K^\pm_i),\quad \sigma(\cZ^\pm_{i,m})=\cZ^\pm_{i,m}(\hV)\quad \forall i\in I, m\in \Z.
\en
\end{prop}
Then we  
 define a functor $\Omega: \gC_{k,l}\to \gD_{k,l}$ by
\be
&&\Omega(\hV,\hpi)=(\Omega_{\hV},\sigma),\quad \Omega({\varphi})=\varphi|_{\Omega_{\hV}}:\Omega_{\hV}\to\Omega_{\hV'}.
\en

\subsection{The functor $\Lambda$}
We define a reverse functor $\Lambda:\gD_{k,l}\to \gC_{k,l}$ as follows. 
Let $(\cW,\sigma)\in \gD_{k,l}$ be a ${\cZ}_{k,l}$-module. We define $U(\gH)$-module 
$\Ind\, \cW$ by 
\be
&&\al_{i,m}\cdot \cW=0,\\
&&\Ind \, \cW=U(\gH)\otimes_{U(\gH_+)}\cW.
\en
Let $\F^{\al}_{k,l}$ be the level-$(k,l)$ Fock module defined in Sec.\ref{Heisenberg}.  We have a natural isomorphism $ 
\F^{\al}_{k,l}\otimes_\C \cW \cong\Ind\ \cW$.  

Identifying $\Ind\, \cW $ with $\F^{\al}_{k,l}\otimes_\C \cW$, 
we define $x^{+'}_j(z), x^{-'}_j(z)\in \cD_{H, \Ind\, \cW}[[z,z^{-1}]]$ by
\be
&&x^{+'}_j(z):=E^-(\al_j,z)^{-1}E^+(\al_j,z)^{-1}\otimes \sigma(\cZ^{+}_{j}(z)) ,\\
&&x^{-'}_j(z):=E^-(\al'_j,z)^{-1}E^+(\al'_j,z)^{-1}\otimes \sigma(\cZ^{-}_{j}(z)).
\en
By a similar argument to the proof of Theorem \ref{zz} one can show that $x^{+'}_j(z)$ and $x^{-'}_j(z)$ gives a level-$(k,l)$ representation of $U_{q,\kappa,p}(\gt)$. 
We hence extend $\hpi: U(\gH)\to \cD_{H,\Ind\, \cW}$ to 
$\hpi: U_{q,\kappa,p}(\g_{tor})\to \cD_{H,\Ind\, \cW}$ as an $H$-algebra homomorphism by 
\be
&&\hpi(x^{+}_j(z))=x^{+'}_j(z),\quad  \hpi(x^{-}_j(z))=x^{-'}_j(z),\\
&&\hpi(q^d)=q^d\otimes 1+1\otimes q^{\sigma(d)}. 
\en
By construction, the latter map is uniquely determined. 
\begin{prop}\lb{IndW}
For $(\cW,\sigma)\in \gD_{k,l}$, there is a unique level-$(k,l)$ $U_{q,\kappa,p}(\gt)$-module 
$(\Ind\, \cW, \hpi)\in \gC_{k,l}$.  
\end{prop}

\begin{dfn}
We define a functor $\Lambda:\gD_{k,l}\to \gC_{k,l}$ by
\begin{itemize}
\item[(i)]\ $\Lambda(\cW,\sigma)=(\Ind\, \cW,\hpi)$
\item[(ii)]\ For a morphism $\varphi:\cW\to \cW'$ in $\gD_{k,l}$, define $\Lambda(\varphi):\Ind\, \cW\to \Ind\, \cW'$ 
to be the induced $U(\gH)$-module map. Then $\Lambda(\varphi)$ is a $U_{q,\kappa,p}(\gt)$-module map. 
\end{itemize}
\end{dfn}

We then obtain the following theorem analogously to the case of the elliptic quantum group $U_{q,p}(\gh)$\cite{FKO}. 
\begin{thm}\lb{IrrIndW}
For $k\in\C^{\times},\ l\in \C$, the two categories $\gC_{k,l}$ and $\gD_{k,l}$ are equivalent by 
the functors $\Omega:\gC_{k,l}\to \gD_k$ and $\Lambda:\gD_{k,l}\to \gC_{k,l}$. 
In particular, the level-$(k,l)$ $U_{q,\kappa,p}(\g_{tor})$-module $\Ind\, \cW=\F^{\al}_{k,l}\otimes_\C \cW \in \gC_{k,l}$ is irreducible 
if and only if $\cW\in \gD_{k,l}$ is an irreducible ${\cZ}_{k,l}$-module.   
\end{thm}

\subsection{The quantum $Z$-algebra associated with $U_{q,\kappa}(\gt)$
}

From the dynamical quantum $Z$-algebra $\cZ^\pm_{i,m}$,  let us define 
$Z^\pm_{i,m}$ by 
\be
&&{Z^+_{j,m}}={\cZ^+_{j,m}}\otimes e^{Q_{j}}, \quad 
{Z^-_{j,m}}={\cZ^-_{j,m}}\otimes 1 
\en
and set $Z^\pm_i(z)=\sum_{m\in \Z}Z^\pm_{i,m}z^{-m}$. 
We also relate $K_i^{\pm}$ to  $k^{\pm1}_i$ in $U_{q,\kappa}(\g_{tor})$ by
\be
&&K_i^{\pm}=k^{\pm1}_i\otimes  e^{-Q_{i}}.
\en
Then one can show that $Z^\pm_{i}(z)$, $k_i^{\pm1}$ satisfy the relations  in  Theorem \ref{zz} except for \eqref{gPZ}-\eqref{gPK} by replacing $\cZ_i^\pm(z,\cV)$ and $K_i^{\pm1}$ with $Z^\pm_i(z)$ and $k_i^{\mp1}$, respectively. 
Namely, the topological algebra generated by $Z^\pm_i(z)$ and $k_i^{\pm 1}$ is a non-dynamical version of the quantum $Z$-algebra.  The latter gives the quantum $Z$-algebra  associated with $U_{q,\kappa}(\gt)$.  Namely one can show the following. 

\begin{prop}\lb{cZHga}
Let us set $a_{j,m}=[d_j]a^\vee_{j,m}\ (j\in I, m\in \Z\backslash\{0\})$. Let $E^\pm(a_{j},z)$ be operators in $U_{q,\kappa}(\gt)[[z,z^{-1}]]$ obtained from $E^\pm(\al_{j},z)$ by replacing $\al_{j,m}$ with $a_{j,m}$. 
We set 
\be
&&X^{+'}_j(z):=E^-(a_j,z)^{-1}E^+(a_j,z)^{-1}
\otimes Z^{+}_{j}(z),\\
&&X^{-'}_j(z):=E^-(a_j,z)E^+(a_j,z)\otimes Z^{-}_{j}(z).
\en
Then $X^{\pm'}_j(z)$, $k_j^{\pm 1}$ generate the quantum toroidal algebra $U_{q,\kappa}(\gt)$.   
\end{prop}

We then obtain the following statement similar to Theorem \ref{IrrIndW}. 
\begin{thm}\lb{IrrW}
A representation  of $U_{q,\kappa}(\gt)$ 
 is irreducible  if and only if the corresponding $Z$-module is irreducible.  
\end{thm}

\subsection{Example
}\lb{level1Ex}

We give an example of the level-$(1,l)$ irreducible induced representation of $U_{q,\kappa,p}(\g_{tor})$  for a simply laced $\g$.

Let $\C[\cQ]$ be the group algebra of the root lattice $\cQ=\oplus_i \Z\al_i$ 
with the central extension:  
\bea
&&e^{ \alpha_i}e^{\pm \alpha_j}=(-1)^{ a_{ij}}\kappa^{\mp m_{ij}}e^{ \pm\alpha_j}e^{ \alpha_i}\qquad 
(i,j\in I). \lb{groupalg}
\ena
Let us consider the fundamental weight $\Lambda_a=\Lambda_0+\bar{\Lambda}_a$ of $\gh$ with $0\leq a\leq N$ for $A^{(1)}_N$, $a=0,1,N-1,N$  for $D_N^{(1)}$, $a=0,1,2$ for $E_6^{(1)}$, $a=0,1$ for $E_7^{(1)}$, $a=0$ for $E_8^{(1)}$. 
\begin{thm}
For $l\in \C$ and  generic  $\mu\in \h^*$,  an  inequivalent set of the irreducible $\cZ_{1,l}$-{modules} is given by 
$\cW(\Lambda_a,\mu)=\FF\otimes_\C e^{\bar{\Lambda}_a}\C[\cQ]\otimes e^{Q_{{\mu}}}\C[\cR_\cQ]$, on which 
  the actions of $\cZ_j^\pm(z)$ 
 are given by 
\bea
&&\cZ_{j}^{+}(z)=e^{ \alpha_j}z^{h_{j}+1}\otimes e^{-Q_j}, \qquad \cZ_{j}^{-}(z)=e^{ -\alpha_j}z^{-h_{j}+1},\lb{qZop}\\
&&K^\pm_j=q^{\pm h_j}\otimes e^{-Q_j}, 
\ena
with  
\be
&&q^{c/2}\cdot 1=q^{1/2}1,\qquad \prod_{j\in I}(K_j^+)^{a_j^\vee} \cdot 1=q^{-l} 1,\\
&&z^{\pm h_{i}}\cdot e^{\alpha_j}e^{\bar\Lambda_a}1=z^{\pm\langle \alpha_j+\bar{\Lambda}_a,h_i\rangle}e^{ \alpha_j}
e^{\bar\Lambda_a} 1 \qquad (i,j\in I).
\en
\end{thm}
\noindent
{\it Proof.}\ One can check the relations in Theorem \ref{zz} with $k=1$  by using 
the following formulas.  
\be
&&\cZ_i^+(z)Z^{\pm}_j(w)=z^{\pm a_{ij}}:\cZ_i^+(z)\cZ^{\pm}_j(w):,\\
&&\cZ^{\pm}_j(w)Z_i^+(z)=w^{\pm a_{ij}}:\cZ^{\pm}_j(w)\cZ_i^+(z):,\\
&&:\cZ_i^+(z)Z^{\pm}_j(w):=(-1)^{a_{ij}}\kappa^{\mp m_{ij}}:\cZ^{\pm}_j(w)\cZ_i^+(z):, 
\en
where
\be
&&:\cZ_i^+(z)\cZ^{+}_j(w):=e^{\al_i}e^{\al_j}z^{h_i+1}w^{ h_j+1}\otimes e^{-Q_i-Q_j}.
\en
etc.. 
In particular, the Serre relations follows from the identity
\be
&&\sum_{\s\in \gS_2}{\mathrm sgn}\s \left(z_{\s(1)}-q^{-1}z_{\s(2)}\right)
\sum_{r=0}^2(-1)^r\left[\mmatrix{2\cr r\cr}\right]\prod_{s=1}^r(w-q^{-1}\kappa^{m_{ij}}z_{\s(s)})\prod_{s=r+1}^2(q^{-1}w-\kappa^{m_{ij}}z_{\s(s)})
=0. 
\en
Note also 
\be
&&{g}(h,P)e^{\sum_{j\in I}a_j^\vee Q_{j}}=e^{\sum_{j\in I}a_j^\vee Q_{j}}{g}(h,P), \quad \forall {g}(h,P)\in \FF. 
\en
\qed

This Theorem is a toroidal and dynamical analogue of the one in \cite{FJ,Jing00}. 

From Theorem \ref{IrrIndW}, we obtain the level-$(1,l)$ irreducible highest weight representation of $U_{q,\kappa,p}(\gt)$ with the highest weight $(\Lambda_a+\mu,\mu)$ by
\bea
&&x^+_i(z)=\exp\left\{\sum_{m>0}\frac{\al_{i,-m}}{[m]}z^m\right\}
\exp\left\{-\sum_{m>0}\frac{\al_{i,m}}{[m]}z^{-m}\right\}\ e^{\al_i}e^{-Q_i} z^{h_i+1}, 
\\
&&x^-_i(z)=\exp\left\{-\sum_{m>0}\frac{\al'_{i,-m}}{[m]}z^m\right\}
\exp\left\{\sum_{m>0}\frac{\al'_{i,m}}{[m]}z^{-m}\right\}\ e^{-\al_i}z^{-h_i+1}
\ena
on 
\bea
&&\FF\otimes_\C\F^\al_{1,l}\otimes \cW(\Lambda_a,\mu).
\ena
The highest weight vector is given by
\be
&&1\otimes e^{\bar{\Lambda}_a}\otimes e^{Q_\mu}.
\en
In the trigonometric case $U_{q,\kappa}(\gl_{N,tor})$, a similar representation was obtained in \cite{Saito}, 
which uses the  extra Cartan elements $H_{i,0}$ in addition to $\partial_{\bar{\al}_i}$ corresponding to our $h_i$ instead of using the  central extension \eqref{groupalg}.

\section{Level-(0,1) Representation of $U_{q,\kappa,p}(\gl_{N,tor})$}\lb{Seclevelzero}
In this section, we restrict ourself to the case 
 $\g=\gln$, \ $N\geq 3$ and consider the level-(0,1) representation i.e. an elliptic and dynamical analogue of the $q$-Fock representation\cite{VV98, STU,FJMM2}. Note that in this representation one has $p=p^*$. 
Hereafter we set $I=\{0,1,\cdots,N-1\}$, 
\be
&&\theta(z)=\theta_p(z), 
\en
and
\be
&&q_1=\kappa q^{-1},\quad q_2=q^2,\quad q_3=\kappa^{-1} q^{-1}. 
\en

\subsection{The vector representation}
For $k\in I , u\in \C^{\times}$, let $V^{(k)}(u)$ be a vector space over $\cM_{H^*}$ with basis $[u]^{(k)}_j, j\in \Z$. 
$V^{(k)}(u)$ has the $\Z^N$-grading defined by
\be
&&\mbox{deg}[u]^{(k)}_j=m (1_0+\cdots+1_{N-1})+1_k+1_{k-1}+\cdots+1_{k-r}
\en
if $j=mN+r$ and $r\in \{0,1,\cdots,N-1\}$. Here $1_i$ denote the standard generators of $\Z^N$ with $1_{i+N}=1_{i}$.  

Define the action of $\phi^\pm_i(z), x^+_i(z), x^-_i(z)\ i\in I $ by
\bea
&&\phi^\pm_i(z)[u]^{(k)}_j=\left\{\mmatrix{q\frac{\theta(q_1^{j+1}q_3u/z)}{\theta(q_1^{j}u/z)}\biggl|_{\pm}e^{-Q_i}[u]^{(k)}_j,\qquad& i+j\equiv k\cr
q^{-1}\frac{\theta(q_1^{j}q_3^{-1}u/z)}{\theta(q_1^{j+1}u/z)}\biggl|_{\pm}e^{-Q_i}[u]^{(k)}_j,\qquad& i+j+1\equiv k\cr
e^{-Q_i}[u]^{(k)}_j\qquad& \mbox{otherwise}\cr}\right.\nn\\
&&x^+_i(z)[u]^{(k)}_j=\left\{\mmatrix{C_+\delta(q_1^{j+1}u/z)e^{-Q_i}[u]^{(k)}_{j+1},\qquad& i+j+1\equiv k\cr
0,\qquad& i+j+1\not\equiv k\cr}\right.\lb{vecrep}\\
&&x^-_i(z)[u]^{(k)}_j=\left\{\mmatrix{C_-\delta(q_1^{j}u/z)[u]^{(k)}_{j-1},\qquad& i+j\equiv k\cr
0,\qquad& i+j\not\equiv k\cr}\right.,\nn
\ena
where 
\bea
&&C_\pm=\frac{(pq^{\pm 2};p)_\infty}{(p;p)_\infty}.
\lb{def:Cpm}
\ena

\begin{lem}
With the action in \eqref{vecrep}, $V^{(k)}(u)$ is an irreducible, tame, $\Z^N$-graded $U_{q,\kappa,p}(\gl_{N,tor})$-module of level (0,0).    We call this the vector representation. 
\end{lem}

\subsection{The $q$-Fock representation}
Let $\la=(\la_1\geq \la_2\geq \cdots)$ be a partition, where $\la_i\geq 0\ (i=1,2,\cdots)$, and   $l(\la)<+\infty$ be its length. 
For $j\in \Z_{\geq 1}$, let $\la\pm \bold{1}_j=(\la_1,\la_2,\cdots,\la_j\pm1,\cdots)$. 
 We identify $\la$ with a Young diagram $\la=\{(x,y)\ |\ 1\leq x\leq l(\la), \ 1\leq y\leq \la_{x}\ \}$. 
Fix $k\in {I}$. For a box $(x,y)\in \la$, we assign  the content $c((x,y))=x-y+k$\footnote{This is called color in \cite{FJMM2}.} modulo $N$. We call a box $(1,1)$ the root of $\la$ and $k$ the root content.   

Let $\mathfrak{F}^{(k)}(u)$ be the subspace
\be
&&\mathfrak{F}^{(k)}(u)\subset V^{(k)}(u)\tot V^{(k)}(uq_2^{-1})\tot V^{(k)}(uq_2^{-2})\tot \cdots
\en
 spanned by vectors 
\be
\ket{\la}^{(k)}_u=[u]^{(k)}_{\la_1-1}\tot [uq_2^{-1}]^{(k)}_{\la_2-2}\tot [uq_2^{-2}]^{(k)}_{\la_3-3}\tot \cdots.
\en 
By repeated application of the comultiplication $\Delta^{op}$  in \eqref{opco1}-\eqref{opco3} to the vector representation, 
we obtain the following level-$(0,1)$ action of $\cU_{q,\kappa,p}(\gl_{N,tor})$ on $\mathfrak{F}^{(k)}(u)$ inductively
\cite{FFJMM}.   
\begin{thm}\lb{zorep}
For $j\in {I}$, the following gives a level-$(0,1)$ irreducible lowest weight representation of $U_{q,\kappa,p}(\gln)$ on $\mathfrak{F}^{(k)}(u)$.
\bea
&&\phi_j^\pm(z)\ket{\la}^{(k)}_u=\prod_{s=1\atop \la_s+j\equiv s+k}^{l(\la)} \left. q\frac{\theta(q_3u_s/z)}{\theta(q_1^{-1}u_s/z)}\right|_\pm\prod_{s=1\atop \la_s+j+1\equiv s+k}^{l(\la)} \left. q^{-1}\frac{\theta(q_1^{-1}q_3^{-1}u_s/z)}{\theta(u_s/z)}\right|_\pm
e^{-Q_{j}}\ket{\la}^{(k)}_u,\lb{phijpm0}\\
&&x^+_j(z)\ket{\la}^{(k)}_u=C_+\sum_{i\in {I}\atop \la_i+j+1\equiv i+k}\delta(u_i/z)A^+_{\la,i}e^{-Q_{j}}
\ket{\la+\bold{1}_i}^{(k)}_u,\lb{ej0}\\
&&x^-_j(z)\ket{\la}^{(k)}_u=C_-\sum_{i\in {I}\atop \la_i+j\equiv i+k}\delta(q_1^{-1}u_i/z)A^-_{\la,i}
\ket{\la-\bold{1}_i}^{(k)}_u,  \lb{fj0}
\ena
where $\equiv$ denotes congruence modulo $N$, and  we set $u_a=q_1^{\la_a}q_3^{a-1} u$  $(a\in I )$, 
\bea
&&A^+_{\la,i}=\prod_{s=1\atop \la_s+j\equiv s+k}^{i-1} q^{-1}\frac{\theta(q_3^{-1}u_i/u_s)}{\theta(q_1u_i/u_s)}\prod_{s=1\atop \la_s+j+1\equiv s+k}^{i-1} q\frac{\theta(q_1q_3u_i/u_s)}{\theta(u_i/u_s)},\lb{def:Ap}\\
&&A^-_{\la,i}=\prod_{s=i+1\atop \la_s+j\equiv s+k}^{\infty} q \frac{\theta(q_1q_3u_s/u_i)}{\theta(u_s/u_i)}
\prod_{s=i+1\atop \la_s+j+1\equiv s+k}^\infty q^{-1}
\frac{\theta(q_3^{-1}u_s/u_i)}{\theta(q_1u_s/u_i)}.\lb{def:Am}
\ena
In \eqref{phijpm0}, $|_\pm$ denotes the expansion direction 
\bea
\left. \frac{\theta(q_3z)(p;p)_\infty^3}{\theta(q_1^{-1}z)\theta(q^{-2})}\right|_+&=&\sum_{n\in \Z}\frac{1}{1-q^{-2}p^n}(q_3z)^n,
\ena
for $|p|<|q_3 z|<1$, and 
\bea
\left. \frac{\theta(q_3z)(p;p)_\infty^3}{\theta(q_1^{-1}z)\theta(q^{-2})}\right|_-
&=&-\frac{\theta(q_3^{-1}/z)(p;p)_\infty^3}{\theta(q_1/z)\theta(q^{2})}\nn\\
&=&\sum_{n\in \Z}\frac{1}{1-q^{2}p^n}(q_3^{-1}/z)^n,
\ena
for $1<|q_3z|<|p^{-1}|$. 
In \eqref{ej0} and \eqref{fj0}, $\ket{\la\pm\bold{1}_i}^{(k)}_u$ vanishes if ${\la\pm \bold{1}_i}$ is  not a partition, respectively. 

The lowest weight vector is $\ket{\emptyset}^{(k)}_u$ and the lowest weight is 
\be
&&\Phi(z)=\left(\left(q^{-1}\frac{\theta(q^2u/z)}{\theta(u/z)}\right)^{\delta_{j,k}}\right)_{j\in I }.
\en
\end{thm}
A direct proof is given in Appendix \ref{Proof01}.

\medskip
A box $(x,y)$ is called addable if $(x,y)\not\in \la$ and $(x-1,y), (x,y-1)\in \la$. 
A box $(x,y)$ is called removable if $(x,y)\in \la$ and $(x+1,y), (x,y+1)\not\in \la$.
For $j\in I $, let $A^{(k)}_j(\la)$ and $R^{(k)}_j(\la)$ denote the set of all addable and removable boxes 
of content $j$, respectively. For $(x,y),  (x',y')\in \la$, let  $(x,y) > (x',y')$ if $c((x,y))> c((x',y'))$.  Then one can 
rewrite the formulas in Theorem \ref{zorep} as follows.
\begin{cor}\lb{lebel0repbox}
\bea
&&\phi_j^\pm(z)\ket{\la}^{(k)}_u=\prod_{R\in R^{(k)}_j(\la)} \left. q\frac{\theta(u_{R}/z)}{\theta(q^2u_R/z)}\right|_\pm
\prod_{A\in A^{(k)}_j(\la)} \left. q^{-1}\frac{\theta(q^4u_A/z)}{\theta(q^2u_A/z)}\right|_\pm
e^{-Q_{j}}\ket{\la}^{(k)}_u,\lb{phijpm0box}\\
&&x^+_j(z)\ket{\la}^{(k)}_u=C_+\sum_{X\in A^{(k)}_j(\la)}\delta(q^2u_X/z)A^+_{\la,X}e^{-Q_{j}}
\ket{\la\cup X}^{(k)}_u,\lb{ej0box}\\
&&x^-_j(z)\ket{\la}^{(k)}_u=C_-\sum_{X\in R^{(k)}_j(\la)}\delta(q^2u_X/z)A^-_{\la,X}
\ket{\la\backslash X}^{(k)}_u,  \lb{fj0box}
\ena
where  $u_{X}=q_1^{y}q_3^{x} u$  for $X=(x,y) \in \la$, and  
\bea
A^+_{\la,X}&=&\prod_{R\in R^{(k)}_j(\la) \atop R<X} q^{-1}\frac{\theta(q^2u_X/u_R)}{\theta(u_X/u_R)}
\prod_{A\in A^{(k)}_j(\la) \atop A<X } q\frac{\theta(q^{-2}u_X/u_A)}{\theta(u_X/u_A)},\lb{def:Apbox}\\
&=&\prod_{R\in R^{(k)}_j(\la) \atop R<X} q^{-1}\frac{\theta(q^2u_X/u_R)}{\theta(u_X/u_R)}
\prod_{A\in A^{(k)}_j(\la\cup X) \atop A<X } q\frac{\theta(q^{-2}u_X/u_A)}{\theta(u_X/u_A)},\lb{Apbox}
\ena
\bea
A^-_{\la,X}&=&\prod_{R\in R^{(k)}_j(\la) \atop R> X} q \frac{\theta(q^{-2}u_R/u_X)}{\theta(u_R/u_X)}
\prod_{A\in A^{(k)}_j(\la) \atop A>X} q^{-1}
\frac{\theta(q^2u_A/u_X)}{\theta(u_A/u_X)}.\lb{def:Ambox}\\
&=&\prod_{R\in R^{(k)}_j(\la\backslash X) \atop R> X} q \frac{\theta(q^{-2}u_R/u_X)}{\theta(u_R/u_X)}
\prod_{A\in A^{(k)}_j(\la) \atop A>X} q^{-1}
\frac{\theta(q^2u_A/u_X)}{\theta(u_A/u_X)}.\lb{def:Ambox}
\ena
In particular,
\bea
&&\prod_{j\in I}K^+_j\ket{\la}^{(k)}_u=q^{|R(\la)|-|A(\la)|}\ket{\la}^{(k)}_u=q^{-1}\ket{\la}^{(k)}_u,\lb{kappa0}
\ena
where we set $|R(\la)|=\sum_{j\in I}|R^{(k)}_j(\la)|$ and $|A(\la)|=\sum_{j\in I}|A^{(k)}_j(\la)|$.
\end{cor}
\noindent
{\it Proof.}\ 
The statement follows from the fact that in \eqref{phijpm0}, \eqref{def:Ap}, \eqref{def:Am}, each factor associated with a box $\not\in R^{(k)}_j(\la)$ in the first product is cancelled out by a factor associated with a box $\not\in A^{(k)}_j(\la)$ in the second product and vice versa. 
Note that if $\la_r=\la_{r+1}$ then $(r,\la_r)\not\in R^{(k)}_j(\la)$ and $(r+1,\la_{r+1}+1)\not\in A^{(k)}_j(\la)$.  
The last equality in \eqref{kappa0} follows from the fact that for any Young diagram the number of addable boxes always exceeds 
the one of removable boxes by one.  
\qed

Finally let us discuss a geometric interpretation of the results. 
For a partition $\la$, let us define  $\bold{d}(\la)=(d_0(\la), d_1(\la),\cdots,d_{N-1}(\la))$ by
 $d_i(\la)=\#$ of boxes having the content $i$ mod $N$. 
Consider the  {quiver variety   ${\cal N}_{\bold{v},\bold{w}}$ of type $A^{(1)}_{N-1}$} with 
 $\bold{v}=\bold{d}(\la)$ and framing $\bold{w}=(w_0,w_1,\cdots,w_{N-1})$ with $w_i=\delta_{i,k}$.  Due to \cite{Nagao}, the trigonometric limit $(p\to 0)$ of the level-$(0,1)$ representation with $k=0$ in Collorary \ref{lebel0repbox} coincides with the representation of $U_{q,\kappa}(\gl_{N,tor})$  on the equivariant 
$\mathrm{K}$-theory $\bigoplus_{\bold{v}\in \N^N}{\mathrm K}_T({\cal N}_{\bold{v},\bold{w}})$ 
with $T= \C^\times \times\C^\times \times (\C^\times)^{|\bold{w}|=1} \ni (q_1,q_3, u)$ 
studied in  \cite{VV99, Na01} except for a dependence on the dynamical parameters. 
 In particular, the vector $\ket{\la}^{(k)}_u$  can be identified with the fixed point class in ${\mathrm K}_T({\cal N}_{\bold{v},\bold{w}})$ labeled by the partition $\la$. Note also that the  level-$(0,1)$ representation of $U_{q,\kappa}(\gl_{N,tor})$ 
 on the $q$-Fock module was constructed in \cite{STU,VV98}. 

We thus reach the following conjecture.

\begin{conj} 
Let $\E_T(X)$ denote the equivariant elliptic cohomology of 
the  quiver variety $X$ defined in \cite{AO}. 
The representation  in Theorem \ref{zorep} gives the level-$(0,1)$ irreducible representation of $U_{q,\kappa,p}(\gl_{N,tor})$ on $\bigoplus_{\bold{v}\in \N^N}{\mathrm E}_T({\cal N}_{\bold{v},\bold{w}})$.
\end{conj}

\bigskip
\section*{Acknowledgements}
The authors would like to thank Yoshihisa Saito for useful discussions. 
 HK is supported by the Grant-in-Aid for Scientific Research (C) 20K03507 JSPS, Japan.
KO is supported by the Grant-in-Aid for Scientific Research (C) 21K03191 JSPS, Japan. 

\appendix
\setcounter{equation}{0}

\begin{appendix}

\section{Lemma \ref{lem:sec4-2}}
\lb{lem:sec4-2}
\begin{align}
& [\alpha_{i,-l}, E^+(\alpha_j,z)]=-\frac{[b_{ij}l]}{l}\frac{1-p^{l}}{1-p^{*l}}\kappa^{l\,m_{ij}}q^{-kl}z^{-l}E^+(\alpha_j,z) \quad (l>0), \\
& [\alpha_{i,l}, E^-(\alpha_j,z)]=-\frac{[b_{ij}l]}{l}\frac{1-p^l}{1-p^{*l}}\kappa^{-l\,m_{ij}}q^{-kl}z^l E^-(\alpha_j,z) \quad (l>0), \\
& [\alpha_{i,-l}, E^+(\alpha'_j,z)]=\frac{[b_{ij}l]}{l}\kappa^{l\,m_{ij}}z^{-l}E^+(\alpha'_j,z) \quad (l>0), \\
& [\alpha_{i,l}, E^-(\alpha'_j,z)]=\frac{[b_{ij}l]}{l}\kappa^{-l\,m_{ij}} z^l E^-(\alpha'_j,z) \quad (l>0), \\
& E^+(\alpha_i,z)E^-(\alpha_j,w)=
\frac{(q^{-b_{ij}}\kappa^{-m_{ij}}\frac{w}{z};q^{2k})_{\infty}(p^*q^{-b_{ij}}\kappa^{-m_{ij}}\frac{w}{z};p^*)_{\infty}}
{(q^{b_{ij}}\kappa^{-m_{ij}}\frac{w}{z};q^{2k})_{\infty}(p^*q^{b_{ij}}\kappa^{-m_{ij}}\frac{w}{z};p^*)_{\infty} }
E^-(\alpha_j,w)E^+(\alpha_i,z), \\
& E^+(\alpha'_i,z)E^-(\alpha'_j,w)=
\frac{(q^{-b_{ij}}q^{2k} \kappa^{-m_{ij}}\frac{w}{z};q^{2k})_{\infty}(pq^{b_{ij}} \kappa^{-m_{ij}}\frac{w}{z};p)_{\infty}}
{(q^{b_{ij}}q^{2k} \kappa^{-m_{ij}}\frac{w}{z};q^{2k})_{\infty}(pq^{-b_{ij}} \kappa^{-m_{ij}} \frac{w}{z};p)_{\infty} }
E^-(\alpha'_j,w)E^+(\alpha'_i,z), \\
& E^+(\alpha_i,z)E^-(\alpha'_j,w)=
\frac{(q^{b_{ij}}q^k \kappa^{-m_{ij}} \frac{w}{z};q^{2k})_{\infty}}{(q^{-b_{ij}}q^k \kappa^{-m_{ij}}\frac{w}{z};q^{2k})_{\infty}} 
E^-(\alpha'_j,w)E^+(\alpha_i,z), \\
& E^+(\alpha'_i,z)E^-(\alpha_j,w)=
\frac{(q^{b_{ij}}q^k \kappa^{-m_{ij}}\frac{w}{z};q^{2k})_{\infty}}{(q^{-b_{ij}}q^k \kappa^{-m_{ij}} \frac{w}{z};q^{2k})_{\infty}} 
E^-(\alpha_j,w)E^+(\alpha'_i,z), \\
& E^+(\alpha_i, z)x_j^+(w)=\frac{(q^{b_{ij}}\kappa^{-m_{ij}}\frac{w}{z};q^{2k})_{\infty}(p^*q^{b_{ij}}\kappa^{-m_{ij}}\frac{w}{z};p^*)_{\infty}}
{(q^{-b_{ij}}\kappa^{-m_{ij}}\frac{w}{z};q^{2k})_{\infty}(p^*q^{-b_{ij}}\kappa^{-m_{ij}}\frac{w}{z};p^*)_{\infty} }
x^+(w,p)E^+(\alpha_i,z), \\
& E^-(\alpha_i, z)x_j^+(w)= \frac{(q^{-b_{ij}}\kappa^{m_{ij}}\frac{z}{w};q^{2k})_{\infty}(p^*q^{-b_{ij}}\kappa^{m_{ij}}\frac{z}{w};p^*)_{\infty}}
{(q^{b_{ij}}\kappa^{m_{ij}}\frac{z}{w};q^{2k})_{\infty}(p^*q^{b_{ij}}
\kappa^{m_{ij}}\frac{z}{w};p^*)_{\infty} }
x_j^+(w)E^-(\alpha_i,z), \\
& E^+(\alpha_i, z)x_j^-(w)=\frac{(q^{-b_{ij}}q^k \kappa^{-m_{ij}}\frac{w}{z};q^{2k})_{\infty}}{(q^{b_{ij}}q^k \kappa^{-m_{ij}}\frac{w}{z};q^{2k})_{\infty}}x_j^-(w)E^+(\alpha_i,z), \\
& E^-(\alpha_i, z)x_j^-(w)=\frac{(q^{b_{ij}}q^k \kappa^{m_{ij}}\frac{z}{w};q^{2k})_{\infty}}{(q^{-b_{ij}}q^k \kappa^{m_{ij}} \frac{z}{w};q^{2k})_{\infty}}x_j^-(w)E^-(\alpha_i,z), \\
& E^+(\alpha'_i, z)x_j^+(w)=\frac{(q^{-b_{ij}}q^k \kappa^{-m_{ij}} \frac{w}{z};q^{2k})_{\infty}}{(q^{b_{ij}}q^k \kappa^{-m_{ij}} \frac{w}{z};q^{2k})_{\infty}} x_j^+(w)E^+(\alpha'_i,z), \\
& E^-(\alpha'_i, z)x_j^+(w)=\frac{(q^{b_{ij}}q^k \kappa^{m_{ij}} \frac{z}{w};q^{2k})_{\infty}}{(q^{-b_{ij}}q^k \kappa^{m_{ij}}\frac{z}{w};q^{2k})_{\infty}} x_j^+(w)E^-(\alpha'_i,z), \\
& E^+(\alpha'_i, z)x_j^-(w)=
\frac{(q^{b_{ij}}q^{2k} \kappa^{-m_{ij}}\frac{w}{z};q^{2k})_{\infty}(pq^{-b_{ij}} \kappa^{-m_{ij}}\frac{w}{z};p)_{\infty}}
{(q^{-b_{ij}}q^{2k} \kappa^{-m_{ij}} \frac{w}{z};q^{2k})_{\infty}(pq^{b_{ij}} \kappa^{-m_{ij}} \frac{w}{z};p)_{\infty} }
 x_j^-(w)E^+(\alpha'_i,z), \\
& E^-(\alpha'_i, z)x_j^-(w)=
\frac{(q^{-b_{ij}}q^{2k} \kappa^{m_{ij}}\frac{z}{w};q^{2k})_{\infty}(pq^{b_{ij}}\kappa^{m_{ij}}\frac{z}{w};p)_{\infty}}
{(q^{b_{ij}}q^{2k} \kappa^{m_{ij}} \frac{z}{w};q^{2k})_{\infty}(pq^{-b_{ij}} \kappa^{m_{ij}}\frac{z}{w};p)_{\infty} }
 x_j^-(w)E^-(\alpha'_i,z).
\end{align}

\section{ Proof of Theorem \ref{zorep}.} \lb{Proof01}
We show that the action in Theorem \ref{zorep} satisfies the 
level-(0,1) relations in Definition \ref{defEQTA}. 
Remember that $p^*=p$ in the level-(0,1) representation. 

For $i\in I$, let us set $u_i=q_1^{\la_i}q_3^{i-1}$ as before.

\noindent
$\bullet$\ \eqref{phixp} with $i=j$
:\\
From \eqref{phijpm0} and \eqref{ej0}, we have 
\be
&&\phi^\pm_j(z)x^+_j(w)\phi^\pm_j(z)^{-1}\ket{\la}^{(k)}_u\\
&&=C_+
\sum_{i\in I\atop \la_i+j+1\equiv i+k}
\prod_{s=1\atop \la_s+j\equiv s+k}^\infty \left. q^{-1}\frac{\theta(q_1^{-1}{u}_s/z)}{\theta(q_3{u}_s/z)}\right|_\pm
\prod_{s=1\atop \la_s+j+1\equiv s+k}^\infty \left. q\frac{\theta({u}_s/z)}{\theta(q_1^{-1}q_3^{-1}{u}_s/z)}\right|_\pm\\
&&\quad\times \delta(u_i/w)A^+_{\la,i}
\prod_{s=1\atop (\la+\bold{1}_i)_s+j\equiv s+k}^\infty \left. q\frac{\theta(q_3\widetilde{u}_s/z)}{\theta(q_1^{-1}\widetilde{u}_s/z)}\right|_\pm
\prod_{s=1\atop (\la+\bold{1}_i)_s+j+1\equiv s+k}^\infty \left. q^{-1}\frac{\theta(q_1^{-1}q_3^{-1}\widetilde{u}_s/z)}{\theta(\widetilde{u}_s/z)}\right|_\pm
e^{-Q_j}\ket{\la+\bold{1}_i}^{(k)}_u,
\en
where we set $\widetilde{u}_s=q_1^{(\la+\bold{1}_i)_s}q_3^{s-1}u=q_1^{\la_s+\delta_{s,i}}q_3^{s-1}u$. 
In side of the summation there exist two $s=i$ factors satisfying $\la_i+j+1\equiv i+k$, i.e.  $q\frac{\theta(u_i/z)}{\theta(q_1^{-1}q_3^{-1}u_i/z)}$  
and $q\frac{\theta(q_1q_3u_i/z)}{\theta(u_i/z)}$ from the products $\displaystyle{\prod_{s=1\atop \la_s+j+1\equiv s+k}^\infty}$ and $\displaystyle{\prod_{s=1\atop (\la+\bold{1}_i)_s+j\equiv s+k}^\infty}$, respectively.  
Hence  we obtain  
\be
&&\phi^\pm_j(z)x^+_j(w)\phi^\pm_j(z)^{-1}\ket{\la}^{(k)}_u\\
&&=q^2\frac{\theta(q_1q_3w/z)}{\theta(q_1^{-1}q_3^{-1}w/z)}C_+
\\
&&\quad\times\sum_{i\in I\atop \la_i+j+1\equiv i+k}
\prod_{s=1\atop {\not=i\atop \la_s+j\equiv s+k}}^\infty \left. 
q^{-1}\frac{\theta(q_1^{-1}u_s/z)}{\theta(q_3u_s/z)}\right|_\pm
\prod_{s=1\atop {\not=i \atop \la_s+j+1\equiv s+k}}^\infty \left. q\frac{\theta(u_s/z)}{\theta(q_1^{-1}q_3^{-1}u_s/z)}\right|_\pm\\
&&\quad\times \delta(u_i/w)A^+_{\la,i}
\prod_{s=1\atop {\not=i\atop \la_s+j\equiv s+k}}^\infty \left. q\frac{\theta(q_3u_s/z)}{\theta(q_1^{-1}u_s/z)}\right|_\pm
\prod_{s=1\atop {\not=i \atop \la_s+j+1\equiv s+k}}^\infty \left. q^{-1}\frac{\theta(q_1^{-1}q_3^{-1}u_s/z)}{\theta(u_s/z)}\right|_\pm
e^{-Q_j}\ket{\la+\bold{1}_i}^{(k)}_u\\
&&=q^2\frac{\theta(q^{-2}w/z)}{\theta(q^2w/z)}x^+_j(w)\ket{\la}^{(k)}_u.
\en

\noindent
$\bullet$\ \eqref{phixp} with $i=j+1$
:
\be
&&\phi^\pm_{j+1}(z)x^+_j(w)\phi^+_{j+1}(z)^{-1}\ket{\la}^{(k)}_u\\
&&=C_+
\sum_{i\in I\atop \la_i+j+1\equiv i+k}
\prod_{s=1\atop \la_s+j+1\equiv s+k}^\infty \left. q^{-1}\frac{\theta(q_1^{-1}{u}_s/z)}{\theta(q_3{u}_s/z)}\right|_\pm
\prod_{s=1\atop \la_s+j+2\equiv s+k}^\infty \left. q\frac{\theta({u}_s/z)}{\theta(q_1^{-1}q_3^{-1}{u}_s/z)}\right|_\pm\\
&&\quad \times \delta(u_i/w)A^+_{\la,i}
\prod_{s=1\atop (\la+\bold{1}_i)_s+j+1\equiv s+k}^\infty \left. q\frac{\theta(q_3\widetilde{u}_s/z)}{\theta(q_1^{-1}\widetilde{u}_s/z)}\right|_\pm
\prod_{s=1\atop (\la+\bold{1}_i)_s+j+2\equiv s+k}^\infty \left. q^{-1}\frac{\theta(q_1^{-1}q_3^{-1}\widetilde{u}_s/z)}{\theta(\widetilde{u}_s/z)}\right|_\pm\\
&&\hspace{6cm}\times e^{-Q_j}\ket{\la+\bold{1}_i}^{(k)}_u.
\en
In side of the summation only the first product has a $s=i$ factor $q^{-1}\frac{\theta(q_1^{-1}u_i/z)}{\theta(q_3u_i/z)}$, we obtain
\be
&&\phi^\pm_{j+1}(z)x^+_j(w)\phi^+_{j+1}(z)^{-1}\ket{\la}^{(k)}_u=q^{-1}\frac{\theta(q\kappa^{-1}w/z)}{\theta(q^{-1}\kappa^{-1}w/z)}x^+_j(w)\ket{\la}^{(k)}_u.
\en
\noindent
One can prove 
\eqref{phixp} with $i=j-1$ 
 similarly.     \\

\noindent
$\bullet$\ \eqref{xpxp} with $i=j$
:\\  
From \eqref{ej0}, we have
\be
&&x^+_j(z)x^+_j(w)\ket{\la}^{(k)}_u\\
&&=C_+^2
\left[
\sum_{i,\ip\in I\atop {\la_i+j+1\equiv i+k \atop {\la_\ip+j+1\equiv \ip+k \atop i<\ip} }}\delta(u_i/w)\delta(u_{i'}/z)
\prod_{s=1\atop \la_s+j\equiv s+k}^{i-1} q^{-1}\frac{\theta(q_3^{-1}u_i/u_s)}{\theta(q_1u_i/u_s)}\prod_{s=1\atop \la_s+j+1\equiv s+k}^{i-1} q\frac{\theta(q_1q_3u_i/u_s)}{\theta(u_i/u_s)}\right.\\
&&\qquad\qquad\qquad\qquad\qquad\times \prod_{\spp=1\atop \la_\spp+j\equiv \spp+k}^{\ip-1} q^{-1}\frac{\theta(q_3^{-1}u_\ip/\widetilde{u}_\spp)}{\theta(q_1u_\ip/\widetilde{u}_\spp)}\prod_{\spp=1\atop \la_\spp+j+1\equiv \spp+k}^{\ip-1} q\frac{\theta(q_1q_3u_\ip/\widetilde{u}_\spp)}{\theta(u_\ip/\widetilde{u}_\spp)}
\\
&&\qquad\qquad\qquad+\sum_{i,\ip\in I\atop {\la_i+j+1\equiv i+k \atop {\la_\ip+j+1\equiv \ip+k \atop i>\ip} }}\delta(u_i/w)\delta(u_{i'}/z)
\prod_{s=1\atop \la_s+j\equiv s+k}^{i-1} q^{-1}\frac{\theta(q_3^{-1}u_i/u_s)}{\theta(q_1u_i/u_s)}\prod_{s=1\atop \la_s+j+1\equiv s+k}^{i-1} q\frac{\theta(q_1q_3u_i/u_s)}{\theta(u_i/u_s)}\\
&&\qquad\qquad\qquad\times \left.\prod_{\spp=1\atop \la_\spp+j\equiv \spp+k}^{\ip-1} q^{-1}\frac{\theta(q_3^{-1}u_\ip/\widetilde{u}_\spp)}{\theta(q_1u_\ip/\widetilde{u}_\spp)}\prod_{\spp=1\atop \la_\spp+j+1\equiv \spp+k}^{\ip-1} q\frac{\theta(q_1q_3u_\ip/\widetilde{u}_\spp)}{\theta(u_\ip/\widetilde{u}_\spp)}\right]
e^{-2Q_j}\ket{\la+\bold{1}_i+\bold{1}_\ip}^{(k)}_u,
\en
where we set $\widetilde{u}_\spp=q_1^{\la_\spp+\delta_{\spp,i}}q_3^{\spp-1}u$. 
In the first term we have a  $\spp=i$  factor $q^{-1}\frac{\theta(q_1^{-1}q_3^{-1}u_\ip/u_i)}{\theta(u_\ip/u_i)}$ in the third product whereas in the second term we have a  $s=i'$  factor $q\frac{\theta(q_1q_3u_i/u_\ip)}{\theta(u_i/u_\ip)}$ in the second product. Hence we obtain
 \be
&&x^+_j(z)x^+_j(w)\ket{\la}^{(k)}_u\nn\\
&&\hspace{-0.5cm}=C_+^2
\left[
q^{-1}\frac{\theta(q^2z/w)}{\theta(z/w)}
\sum_{i,\ip\in I\atop {\la_i+j+1\equiv i+k \atop {\la_\ip+j+1\equiv \ip+k \atop i<\ip} }}\delta(u_i/w)\delta(u_{i'}/z)
\prod_{s=1\atop \la_s+j\equiv s+k}^{i-1} q^{-1}\frac{\theta(q_3^{-1}u_i/u_s)}{\theta(q_1u_i/u_s)}\prod_{s=1\atop \la_s+j+1\equiv s+k}^{i-1} q\frac{\theta(q_1q_3u_i/u_s)}{\theta(u_i/u_s)}
\right.
\\
&&\qquad\qquad\qquad\qquad\qquad\times \prod_{\spp=1\atop {\not=i \atop \la_\spp+j\equiv \spp+k}}^{\ip-1} q^{-1}\frac{\theta(q_3^{-1}u_\ip/\widetilde{u}_\spp)}{\theta(q_1u_\ip/\widetilde{u}_\spp)}\prod_{\spp=1\atop \la_\spp+j+1\equiv \spp+k}^{\ip-1} q\frac{\theta(q_1q_3u_\ip/\widetilde{u}_\spp)}{\theta(u_\ip/\widetilde{u}_\spp)}
\en
\be
&&\qquad+q\frac{\theta(q^{-2}w/z)}{\theta(w/z)}
\sum_{i,\ip\in I\atop {\la_i+j+1\equiv i+k \atop {\la_\ip+j+1\equiv \ip+k \atop i>\ip} }}\delta(u_i/w)\delta(u_{i'}/z)
\prod_{s=1\atop \la_s+j\equiv s+k}^{i-1} q^{-1}\frac{\theta(q_3^{-1}u_i/u_s)}{\theta(q_1u_i/u_s)}\prod_{s=1\atop {\not=\ip \atop  \la_s+j+1\equiv s+k}}^{i-1} q\frac{\theta(q_1q_3u_i/u_s)}{\theta(u_i/u_s)}\\
&&\qquad\qquad\qquad\times \left.\prod_{\spp=1\atop \la_\spp+j\equiv \spp+k}^{\ip-1} q^{-1}\frac{\theta(q_3^{-1}u_\ip/\widetilde{u}_\spp)}{\theta(q_1u_\ip/\widetilde{u}_\spp)}\prod_{\spp=1\atop \la_\spp+j+1\equiv \spp+k}^{\ip-1} q\frac{\theta(q_1q_3u_\ip/\widetilde{u}_\spp)}{\theta(u_\ip/\widetilde{u}_\spp)}\right]
e^{-2Q_j}\ket{\la+\bold{1}_i+\bold{1}_\ip}^{(k)}_u.
\en
From this expression by exchanging $z$ and $w$ as well as $i$ and $\ip$, simultaneously, 
we obtain a similar formula for $x^+_j(w)x^+_j(z)\ket{\la}^{(k)}_u$. Comparing these two,
 one obtains 
\be
&&x^+_j(z)x^+_j(w)\ket{\la}^{(k)}_u=q^2\frac{\theta(q^{-2}w/z)}{\theta(q^2w/z)} x^+_j(w)x^+_j(z)\ket{\la}^{(k)}_u.
\en

\noindent
$\bullet$\ \eqref{xpxp} 
for $x^+_j(z)$ and $x^+_{j+1}(w)$
:\\
From \eqref{ej0}, we have
 \be
&&q^{-1}\frac{\theta(q_1^{-1}z/w)}{\theta(q_3z/w)}x^+_j(z)x^+_{j+1}(w)\ket{\la}^{(k)}_u\\
&&\hspace{-0.5cm}=C_+^2
\left[
q^{-1}\frac{\theta(q_1^{-1}z/w)}{\theta(q_3z/w)}\sum_{i\in I\atop \la_i+j+2\equiv i+k }
\delta(u_i/w)\delta(q_1u_{i}/z)
\prod_{s=1\atop \la_s+j+1\equiv s+k}^{i-1} q^{-1}\frac{\theta(q_3^{-1}u_i/u_s)}{\theta(q_1u_i/u_s)}\prod_{s=1\atop \la_s+j+2\equiv s+k}^{i-1} q\frac{\theta(q_1q_3u_i/u_s)}{\theta(u_i/u_s)}\right.
\\
&&\qquad\qquad\qquad\times \prod_{\spp=1\atop { \la_\spp+j\equiv \spp+k}}^{i-1} q^{-1}\frac{\theta(q_1q_3^{-1}u_i/{u}_\spp)}{\theta(q_1^2u_\ip/{u}_\spp)}\prod_{\spp=1\atop \la_\spp+j+1\equiv \spp+k}^{i-1} q\frac{\theta(q_1^2q_3u_\ip/{u}_\spp)}{\theta(q_1u_\ip/{u}_\spp)}
\ e^{-Q_j-Q_{j+1}}\ket{\la+2\bold{1}_i}^{(k)}_u
\\
&&+
\sum_{i,\ip\in I\atop {\la_i+j+2\equiv i+k \atop {\la_\ip+j+1\equiv \ip+k \atop i<\ip} }}\delta(u_i/w)\delta(u_{i'}/z)
\prod_{s=1\atop \la_s+j+1\equiv s+k}^{i-1} q^{-1}\frac{\theta(q_3^{-1}u_i/u_s)}{\theta(q_1u_i/u_s)}\prod_{s=1\atop \la_s+j+2\equiv s+k}^{i-1} q\frac{\theta(q_1q_3u_i/u_s)}{\theta(u_i/u_s)}
\\
&&\qquad\qquad\qquad\times \prod_{\spp=1\atop { \la_\spp+j\equiv \spp+k}}^{\ip-1} q^{-1}\frac{\theta(q_3^{-1}u_\ip/{u}_\spp)}{\theta(q_1u_\ip/{u}_\spp)}\prod_{\spp=1\atop {\not=i \atop \la_\spp+j+1\equiv \spp+k}}^{\ip-1} q\frac{\theta(q_1q_3u_\ip/{u}_\spp)}{\theta(u_\ip/{u}_\spp)}
\ e^{-Q_j-Q_{j+1}}\ket{\la+\bold{1}_i+\bold{1}_\ip}^{(k)}_u
\\
&&\qquad+
\sum_{i,\ip\in I\atop {\la_i+j+2\equiv i+k \atop {\la_\ip+j+1\equiv \ip+k \atop i>\ip} }}\delta(u_i/w)\delta(u_{i'}/z)
\prod_{s=1\atop {\not=\ip \atop \la_s+j+1\equiv s+k}}^{i-1} q^{-1}\frac{\theta(q_3^{-1}u_i/u_s)}{\theta(q_1u_i/u_s)}\prod_{s=1\atop { 
 \la_s+j+2\equiv s+k}}^{i-1} q\frac{\theta(q_1q_3u_i/u_s)}{\theta(u_i/u_s)}\\
&&\qquad\qquad\qquad\times \left.\prod_{\spp=1\atop \la_\spp+j\equiv \spp+k}^{\ip-1} q^{-1}\frac{\theta(q_3^{-1}u_\ip/{u}_\spp)}{\theta(q_1u_\ip/{u}_\spp)}\prod_{\spp=1\atop \la_\spp+j+1\equiv \spp+k}^{\ip-1} q\frac{\theta(q_1q_3u_\ip/{u}_\spp)}{\theta(u_\ip/{u}_\spp)}\ 
e^{-Q_j-Q_{j+1}}\ket{\la+\bold{1}_i+\bold{1}_\ip}^{(k)}_u\right].
\en
Then the first term vanishes because of $\theta(q_1^{-1}z/w)\delta(u_i/w)\delta(q_1u_{i}/z)=0$. 
We hence obtain
\be
&&q^{-1}\frac{\theta(q_1^{-1}z/w)}{\theta(q_3z/w)}x^+_j(z)x^+_{j+1}(w)\ket{\la}^{(k)}_u\\
&&\hspace{-0.5cm}=C_+^2
\left[
\sum_{i,\ip\in I\atop {\la_i+j+2\equiv i+k \atop {\la_\ip+j+1\equiv \ip+k \atop i<\ip} }}\delta(u_i/w)\delta(u_{i'}/z)
\prod_{s=1\atop \la_s+j+1\equiv s+k}^{i-1} q^{-1}\frac{\theta(q_3^{-1}u_i/u_s)}{\theta(q_1u_i/u_s)}\prod_{s=1\atop \la_s+j+2\equiv s+k}^{i-1} q\frac{\theta(q_1q_3u_i/u_s)}{\theta(u_i/u_s)}\right.
\\
&&\qquad\qquad\qquad\times \prod_{\spp=1\atop { \la_\spp+j\equiv \spp+k}}^{\ip-1} q^{-1}\frac{\theta(q_3^{-1}u_\ip/{u}_\spp)}{\theta(q_1u_\ip/{u}_\spp)}\prod_{\spp=1\atop {\not=i \atop \la_\spp+j+1\equiv \spp+k}}^{\ip-1} q\frac{\theta(q_1q_3u_\ip/{u}_\spp)}{\theta(u_\ip/{u}_\spp)}
\ e^{-Q_j-Q_{j+1}}\ket{\la+\bold{1}_i+\bold{1}_\ip}^{(k)}_u
\\
&&\qquad+
\sum_{i,\ip\in I\atop {\la_i+j+2\equiv i+k \atop {\la_\ip+j+1\equiv \ip+k \atop i>\ip} }}\delta(u_i/w)\delta(u_{i'}/z)
\prod_{s=1\atop {\not=\ip \atop \la_s+j+1\equiv s+k}}^{i-1} q^{-1}\frac{\theta(q_3^{-1}u_i/u_s)}{\theta(q_1u_i/u_s)}\prod_{s=1\atop { 
 \la_s+j+2\equiv s+k}}^{i-1} q\frac{\theta(q_1q_3u_i/u_s)}{\theta(u_i/u_s)}\\
&&\qquad\qquad\qquad\times \left.\prod_{\spp=1\atop \la_\spp+j\equiv \spp+k}^{\ip-1} q^{-1}\frac{\theta(q_3^{-1}u_\ip/{u}_\spp)}{\theta(q_1u_\ip/{u}_\spp)}\prod_{\spp=1\atop \la_\spp+j+1\equiv \spp+k}^{\ip-1} q\frac{\theta(q_1q_3u_\ip/{u}_\spp)}{\theta(u_\ip/{u}_\spp)}\ 
e^{-Q_j-Q_{j+1}}\ket{\la+\bold{1}_i+\bold{1}_\ip}^{(k)}_u\right]\\
&&=x^+_{j+1}(w)x^+_j(z)\ket{\la}^{(k)}_u.
\en
One can  prove \eqref{xpxp} 
for $x^+_j(z)$ and $x^+_{j-1}(w)$ 
 similarly. 

\noindent
$\bullet$\ \eqref{xpxm} 
:\\
From \eqref{fj0} and \eqref{ej0}, we have
\be
&&x^+_j(z)x^-_j(w)\ket{\la}^{(k)}_u\\
&&=C_+C_-\sum_{i\in I\atop \la_i+j\equiv i+k}\sum_{\ip\in I\atop (\la-\bold{1}_i)_\ip+j+1\equiv \ip+k}\delta(q_1^{-1}u_i/w)\delta(q_1^{-\delta_{i,\ip}}{u}_\ip/z)A^-_{\la,i}A^+_{\la-\bold{1}_i,\ip}e^{-Q_j}\ket{\la-\bold{1}_i+\bold{1}_\ip}^{(k)}_u\\
&&=C_+C_-\sum_{i\in I\atop \la_i+j\equiv i+k}\delta(q_1^{-1}u_i/w)\delta(q_1^{-1}{u}_i/z)A^-_{\la,i}A^+_{\la-\bold{1}_i,i}e^{-Q_j}\ket{\la}^{(k)}_u\\
&&\quad+CC'\sum_{i,\ip\in I\atop {i\not=\ip \atop {\la_i+j\equiv i+k
\atop \la_\ip+j+1\equiv \ip+k
}}
}\delta(q_1^{-1}u_i/w)\delta({u}_\ip/z)A^-_{\la,i}A^+_{\la-\bold{1}_i,\ip}e^{-Q_j}\ket{\la-\bold{1}_i+\bold{1}_\ip}^{(k)}_u.
\en
Note that 
\bea
C_+C_-&=&\frac{q}{q-q^{-1}}\frac{\theta(q^{-2})}{(p;p)_\infty^3},\lb{CCp}
\ena
\bea
A^+_{\la-\bold{1}_i,i}&=&\prod_{s=1\atop { \la_s+j\equiv s+k}}^{i-1} q^{-1}\frac{\theta(q_1^{-1}q_3^{-1}u_i/u_s)}{\theta(u_i/u_s)}\prod_{s=1\atop { 
 \la_s+j+1\equiv s+k}}^{i-1} q\frac{\theta(q_3u_i/u_s)}{\theta(q_1^{-1}u_i/u_s)},\lb{Apii}
\\
A^+_{\la-\bold{1}_i,\ip}&=&\prod_{s=1\atop { \la_s+j\equiv s+k}}^{\ip-1} q^{-1}\frac{\theta(q_3^{-1}u_\ip/u_s)}{\theta(q_1u_\ip/u_s)}\prod_{s=1\atop { 
 \la_s+j+1\equiv s+k}}^{i'-1} q\frac{\theta(q_1q_3u_\ip/u_s)}{\theta(u_\ip/u_s)}\ \qquad\quad\qquad\qquad\qquad \mbox{for } \ i>\ip\nn\\
  &=&\prod_{s=1\atop { \la_s+j\equiv s+k}}^{\ip-1} q^{-1}\frac{\theta(q_3^{-1}u_\ip/u_s)}{\theta(q_1u_\ip/u_s)}
\prod_{s=1\atop { \not=i \atop 
 \la_s+j+1\equiv s+k}
 }^{\ip-1} q\frac{\theta(q_1q_3u_\ip/u_s)}{\theta(u_\ip/u_s)}\times q\frac{\theta(q_1^2q_3u_\ip/u_i)}{\theta(q_1u_\ip/u_i)}
  \qquad \mbox{for } \ i<\ip \nn\\
  &&\lb{Apiip}
\ena
Similarly, one has
\be
&&x^-_j(w)x^+_j(z)\ket{\la}^{(k)}_u\\
&&=C_+C_-\sum_{\ip\in I\atop \la_\ip+j+1\equiv \ip+k}\sum_{i\in I\atop (\la+\bold{1}_\ip)_i+j\equiv i+k}\delta(q_1^{\delta_{i,\ip}-1}u_i/w)\delta({u}_\ip/z)A^+_{\la,i'}A^-_{\la+\bold{1}_\ip,i}e^{-Q_j}\ket{\la-\bold{1}_i+\bold{1}_\ip}^{(k)}_u\\
&&=C_+C_-\sum_{i\in I\atop \la_i+j+1\equiv i+k}\delta(u_i/w)\delta({u}_i/z)A^+_{\la,i}A^+_{\la+\bold{1}_i,i}e^{-Q_j}\ket{\la}^{(k)}_u\\
&&\quad +C_+C_-\sum_{i,\ip\in I\atop {i\not=\ip \atop {\la_i+j\equiv i+k
\atop \la_\ip+j+1\equiv \ip+k
}}
}\delta(q_1^{-1}u_i/w)\delta({u}_\ip/z)A^+_{\la,i'}A^-_{\la+\bold{1}_\ip,i}e^{-Q_j}\ket{\la-\bold{1}_i+\bold{1}_\ip}^{(k)}_u,
\en
with 
\bea
 A^-_{\la+\bold{1}_i,i}&=&\prod_{s=i+1\atop { \la_s+j\equiv s+k}}^{\infty} q\frac{\theta(q_3u_s/u_i)}{\theta(q_1^{-1}u_s/u_i)}
 \prod_{s=i+1\atop { 
 \la_s+j+1\equiv s+k}}^{\infty} q^{-1}\frac{\theta(q_1^{-1}q_3^{-1}u_s/u_i)}{\theta(u_s/u_i)},\lb{Amii}\\
A^-_{\la+\bold{1}_\ip,i}
&=&\prod_{s=i+1\atop { \la_s+j\equiv s+k}}^{\infty} q\frac{\theta(q_1q_3u_s/u_i)}{\theta(u_s/u_i)}
\prod_{s=i+1\atop { 
 \la_s+j+1\equiv s+k}}^{\infty} q^{-1}\frac{\theta(q_3^{-1}u_s/u_i)}{\theta(q_1u_s/u_i)}\ \qquad\quad\qquad\qquad\qquad \mbox{for } \ i>\ip
 \nn\\
 &=&\prod_{s=i+1 \atop
{\not=\ip  \atop { \la_s+j\equiv s+k}}
}^{\infty} q\frac{\theta(q_1q_3u_s/u_i)}{\theta(u_s/u_i)}
\prod_{s=i+1\atop {  
 \la_s+j+1\equiv s+k}
 }^{\infty} q^{-1}\frac{\theta(q_3^{-1}u_s/u_i)}{\theta(q_1u_s/u_i)}\times q\frac{\theta(q_1^2q_3u_\ip/u_i)}{\theta(q_1u_\ip/u_i)}
  \qquad \mbox{for } \ i<\ip\nn\\
  &&\lb{Amiip}
\ena
Therefore one obtains
\be
&&[x^+_j(z),x^-_j(w)]\ket{\la}^{(k)}_u\\
&&=C_+C_-\delta(z/w)
\left(\sum_{i\in I\atop \la_i+j\equiv i+k}\delta(q_1^{-1}{u}_i/z)A^-_{\la,i}A^+_{\la-\bold{1}_i,i}-  
\sum_{i\in I\atop \la_i+j+1\equiv i+k}\delta(q_1^{-1}{u}_i/z)A^+_{\la,i}A^+_{\la+\bold{1}_i,i}   \right)e^{-Q_j}\ket{\la}^{(k)}_u
\en
\be
&&\quad+C_+C_-\sum_{i,\ip\in I\atop {i\not=\ip \atop {\la_i+j\equiv i+k 
\atop \la_\ip+j+1\equiv \ip+k
}
}
}\delta(q_1^{-1}u_i/w)\delta({u}_\ip/z)\left(
A^-_{\la,i}A^+_{\la-\bold{1}_i,\ip}- A^+_{\la,i'}A^-_{\la+\bold{1}_\ip,i}    \right)
e^{-Q_j}\ket{\la-\bold{1}_i+\bold{1}_\ip}^{(k)}_u. 
\en
From \eqref{def:Ap}, \eqref{def:Am}, \eqref{Apiip} and \eqref{Amiip}, the second term vanishes. 
In fact, for $i>i'$ satisfying $\la_i+j\equiv i+k,   \la_\ip+j+1\equiv \ip+k$, we have
\be
A^-_{\la,i}A^+_{\la-\bold{1}_i,\ip}
&=&\prod_{s=i+1 \atop \la_s+j\equiv s+k}^{l(\la)}q\frac{\theta(q_1q_3u_s/u_i)}{\theta(u_s/u_i)}
\prod_{s=i+1 \atop \la_s+j+1\equiv s+k}^{l(\la)}q^{-1}\frac{\theta(q_3^{-1}u_s/u_i)}{\theta(q_1u_s/u_i)}\\
&&\times\prod_{s'=1\atop { \la_\spp+j\equiv s'+k}}^{\ip-1} q^{-1}\frac{\theta(q_3^{-1}u_\ip/u_\spp)}{\theta(q_1u_\ip/u_\spp)}
\prod_{s'=1\atop { 
 \la_\spp+j+1\equiv s'+k}}^{i'-1} q\frac{\theta(q_1q_3u_\ip/u_\spp)}{\theta(u_\ip/u_\spp)},
 \en
\be
A^+_{\la,i'}A^-_{\la+\bold{1}_\ip,i}
&=&\prod_{s'=1 \atop \la_\spp+j\equiv s'+k}^{i'-1}q^{-1}\frac{\theta(q_3^{-1}u_\ip/u_\spp)}{\theta(q_1u_\ip/u_\spp)}
\prod_{s'=1 \atop \la_\spp+j+1\equiv s'+k}^{i'-1}q\frac{\theta(q_1q_3u_\ip/u_\spp)}{\theta(u_\ip/u_\spp)}\\
&&\times \prod_{s=i+1\atop { \la_s+j\equiv s+k}}^{l(\la)} q\frac{\theta(q_1q_3u_s/u_i)}{\theta(u_s/u_i)}
\prod_{s=i+1\atop { 
 \la_s+j+1\equiv s+k}}^{l(\la)} q^{-1}\frac{\theta(q_3^{-1}u_s/u_i)}{\theta(q_1u_s/u_i)}, 
\en
whereas, for $i<i'$  satisfying $\la_i+j\equiv i+k,   \la_\ip+j+1\equiv \ip+k$, we have 
\be
A^-_{\la,i}A^+_{\la-\bold{1}_i,\ip}
&=&\prod_{s=i+1 \atop {(\not=i') \atop \la_s+j\equiv s+k}}^{l(\la)}q\frac{\theta(q_1q_3u_s/u_i)}{\theta(u_s/u_i)}
\prod_{s=i+1 \atop \la_s+j+1\equiv s+k}^{l(\la)}q^{-1}\frac{\theta(q_3^{-1}u_s/u_i)}{\theta(q_1u_s/u_i)}\\
&&\times
\prod_{s=1\atop { \la_s+j\equiv s+k}}^{\ip-1} q^{-1}\frac{\theta(q_3^{-1}u_\ip/u_s)}{\theta(q_1u_\ip/u_s)}
\prod_{s=1\atop { \not=i \atop 
 \la_s+j+1\equiv s+k}
 }^{\ip-1} q\frac{\theta(q_1q_3u_\ip/u_s)}{\theta(u_\ip/u_s)}\times q\frac{\theta(q_1^2q_3u_\ip/u_i)}{\theta(q_1u_\ip/u_i)},
\en
\be
A^+_{\la,i'}A^-_{\la+\bold{1}_\ip,i}
&=&\prod_{s'=1 \atop \la_\spp+j\equiv s'+k}^{i'-1}q^{-1}\frac{\theta(q_3^{-1}u_\ip/u_\spp)}{\theta(q_1u_\ip/u_\spp)}
\prod_{s'=1 \atop {(\not=i) \atop \la_\spp+j+1\equiv s'+k}}^{i'-1}q\frac{\theta(q_1q_3u_\ip/u_\spp)}{\theta(u_\ip/u_\spp)}\\
&&\times \prod_{s=i+1\atop {\not=i' \atop \la_s+j\equiv s+k}}^{l(\la)} q\frac{\theta(q_1q_3u_s/u_i)}{\theta(u_s/u_i)}
\prod_{s=i+1\atop { 
 \la_s+j+1\equiv s+k}}^{l(\la)} q^{-1}\frac{\theta(q_3^{-1}u_s/u_i)}{\theta(q_1u_s/u_i)}\times q\frac{\theta(q_1^2q_3u_\ip/u_i)}{\theta(q_1 u_\ip/u_i)}.
\en
Hence from \eqref{def:Ap}, \eqref{def:Am}, \eqref{CCp}, \eqref{Apii} and \eqref{Amii} we obtain
\bea
&&[x^+_j(z),x^-_j(w)]\ket{\la}^{(k)}_u\nn\\
&&=\frac{q}{q-q^{-1}}\frac{\theta(q^{-2})}{(p;p)_\infty^2}\delta(z/w)\nn\\
&& \times \left(\sum_{i\in I\atop \la_i+j\equiv i+k}\delta(q_1^{-1}{u}_i/z)
 \prod_{s=1\atop { \not=i  \atop 
 \la_s+j\equiv s+k}}^{l(\la)} q^{-1}\frac{\theta(q_1^{-1}q_3^{-1}u_i/u_s)}{\theta(u_i/u_s)}
\prod_{s=1\atop {(\not= i) \atop \la_s+j+1\equiv s+k}}^{l(\la)} q\frac{\theta(q_3u_i/u_s)}{\theta(q_1^{-1}u_i/u_s)}
\right.\nn\\
&&\left.-\sum_{i\in I\atop \la_i+j+1\equiv i+k}\delta({u}_i/z)
\prod_{s=1\atop {(\not= i) \atop \la_s+j\equiv s+k}}^{l(\la)} q^{-1}\frac{\theta(q_3^{-1}u_i/u_s)}{\theta(q_1u_i/u_s)}
 \prod_{s=1\atop { \not=i  \atop 
 \la_s+j+1\equiv s+k}}^{l(\la)} q\frac{\theta(q_1q_3u_i/u_s)}{\theta(u_i/u_s)}\right)e^{-Q_j}\ket{\la}^{(k)}_u.\lb{comejfj}
\ena

Next let us consider \eqref{phijpm0} and use the following partial fraction expansion.  
\begin{lem}\cite{Rosengren}\lb{pfexp}
The following formula  holds. 
\be
\prod_{s=1}^n\frac{\theta(b_s/t)}{\theta(a_s/t)}=\frac{1}{\theta(b_{n+1}/t)}\sum_{i=1}^n\frac{\theta(a_i/b_{n+1})}{\theta(a_i/t)}\frac{\prod_{s=1}^{n}\theta(a_i/b_s)}{\prod_{s=1\atop \not=i}^n\theta(a_i/a_s)}
\en
with $a_1\cdots a_n t=b_1\cdots b_{n+1}$. 
\end{lem}
Applying this formula with taking $t=z$,
\be
&&a_s=q_1^{-1}u_s, \quad b_s=q_3u_s \qquad \mbox{for}\ \la_s+j\equiv s+k \\
&&a_s=u_s, \quad b_s=q_1^{-1}q_3^{-1}u_s \qquad \mbox{for}\ \la_s+j+1\equiv s+k, 
\en
$1\leq s\leq l(\la)$, and
\be
&&b_{n+1}=Q^{-1}z,\qquad Q=(q_1q_3)^{{\sum'}_{s=1}^{l(\la)}-{\sum''}_{s=1}^{l(\la)}}=q^{-2{\sum'}_{s=1}^{l(\la)}+2{\sum''}_{s=1}^{l(\la)}},
\en
where $\sum'_s$ and $\sum''_s$ denote a summation over $s$ satisfying $\la_s+j\equiv s+k$ and $\la_s+j+1\equiv s+k$, respectively,
we obtain 
\be
&&\prod_{s=1\atop \la_s+j\equiv s+k}^{l(\la)}
\left.
 q\frac{\theta(q_3u_s/z)}{\theta(q_1^{-1}u_s/z)}
  \right|_\pm
\prod_{s=1\atop \la_s+j+1\equiv s+k}^{l(\la)}
\left. 
q^{-1}\frac{\theta(q_1^{-1}q_3^{-1}u_s/z)}{\theta(u_s/z)}
\right|_\pm   
\\
&&=-{q^{-{\sum'}_{s=1}^{l(\la)}+{\sum''}_{s=1}^{l(\la)}}}\\
&&\times\left[
\sum_{i=1\atop \la_i+j\equiv i+k}^{l(\la)}\theta(q_1^{-1}q_3^{-1})
\left.\frac{\theta(Qq_1^{-1}u_i/z)}{{\theta(Q)}\theta(q_1^{-1}u_i/z)}\right|_{\pm}
\prod_{s=1\atop {\not=i \atop \la_s+j\equiv s+k}}^{l(\la)}\frac{\theta(q_1^{-1}q_3^{-1}u_i/u_s)}{\theta(u_i/u_s)}
\prod_{s=1\atop {(\not=i) \atop \la_s+j+1\equiv s+k}}^{l(\la)}\frac{\theta(q_3u_i/u_s)}{\theta(q_1^{-1}u_i/u_s)}\right.
\en
\be
&&\left.\qquad+\sum_{i=1\atop \la_i+j+1\equiv i+k}^{l(\la)}\theta(q_1q_3)
\left.\frac{\theta(Qu_i/z)}{{\theta(Q)}\theta(u_i/z)}\right|_{\pm}
\prod_{s=1\atop {(\not=i) \atop \la_s+j\equiv s+k}}^{l(\la)}\frac{\theta(q_3^{-1}u_i/u_s)}{\theta(q_1u_i/u_s)}
\prod_{s=1\atop {\not=i \atop \la_s+j+1\equiv s+k}}^{l(\la)}\frac{\theta(q_1q_3u_i/u_s)}{\theta(u_i/u_s)}\right].\\
\en
Then using the formula
\be
&&\left.\frac{\theta(wz)}{\theta(w)\theta(z)}\right|_+- \left.\frac{\theta(wz)}{\theta(w)\theta(z)}\right|_-=-\frac{1}{(p;p)^2_\infty}\delta(z)
\en
for $w\in \C^{\times}$, 
we obtain 
\be
&&\left(\phi^+_j(z)-\phi^-_j(z)\right)\ket{\la}^{(k)}_u\\
&&=\frac{q \theta(q_1q_3) }{(p;p)^2_\infty}q^{-{\sum'}_{s=1}^{l(\la)}+{\sum''}_{s=1}^{l(\la)}}\\
&&\times\left[
q\sum_{i=1\atop \la_i+j\equiv i+k}^{l(\la)}
\delta(q_1^{-1}u_i/z)
\prod_{s=1\atop {\not=i \atop \la_s+j\equiv s+k}}^{l(\la)}\frac{\theta(q_1^{-1}q_3^{-1}u_i/u_s)}{\theta(u_i/u_s)}
\prod_{s=1\atop {(\not=i) \atop \la_s+j+1\equiv s+k}}^{l(\la)}\frac{\theta(q_3u_i/u_s)}{\theta(q_1^{-1}u_i/u_s)}\right.\\
&&\left.\qquad-q^{-1}\sum_{i=1\atop \la_i+j+1\equiv i+k}^{l(\la)}
\delta(u_i/z)
\prod_{s=1\atop {(\not=i) \atop \la_s+j\equiv s+k}}^{l(\la)}\frac{\theta(q_3^{-1}u_i/u_s)}{\theta(q_1u_i/u_s)}
\prod_{s=1\atop {\not=i \atop \la_s+j+1\equiv s+k}}^{l(\la)}\frac{\theta(q_1q_3u_i/u_s)}{\theta(u_i/u_s)}\right]e^{-Q_j}\ket{\la}^{(k)}_u\\
&&=\frac{ q \theta(q_1q_3)}{(p;p)^2_\infty}
\left[
\sum_{i=1\atop \la_i+j\equiv i+k}^{l(\la)}
\delta(q_1^{-1}u_i/z)
\prod_{s=1\atop {\not=i \atop \la_s+j\equiv s+k}}^{l(\la)}q^{-1}\frac{\theta(q_1^{-1}q_3^{-1}u_i/u_s)}{\theta(u_i/u_s)}
\prod_{s=1\atop {(\not=i) \atop \la_s+j+1\equiv s+k}}^{l(\la)}q\frac{\theta(q_3u_i/u_s)}{\theta(q_1^{-1}u_i/u_s)}\right.\\
&&\left.\qquad-\sum_{i=1\atop \la_i+j+1\equiv i+k}^{l(\la)}
\delta(u_i/z)
\prod_{s=1\atop {(\not=i) \atop \la_s+j\equiv s+k}}^{l(\la)}q^{-1}\frac{\theta(q_3^{-1}u_i/u_s)}{\theta(q_1u_i/u_s)}
\prod_{s=1\atop {\not=i \atop \la_s+j+1\equiv s+k}}^{l(\la)}q\frac{\theta(q_1q_3u_i/u_s)}{\theta(u_i/u_s)}\right]e^{-Q_j}\ket{\la}^{(k)}_u.
\en
Comparing this with \eqref{comejfj}, we obtain
\be
&&[x^+_j(z),x^-_j(w)]\ket{\la}^{(k)}_u=\frac{1}{q-q^{-1}}\delta(z/w)\left(\phi^+_j(z)-\phi^-_j(z)\right)\ket{\la}^{(k)}_u.
\en

\noindent
$\bullet$\ Serre relation \eqref{Serrexp}: \\
Note that \eqref{Serrexp} can be rewritten as 
\bea
&&\frac{(pq^2z_2/z_1;p)_\infty}{(pq^{-2}z_2/z_1;p)_\infty}\prod_{a=1,2}\frac{(pq^{-1}\kappa^{\pm1}w/z_a;p)_\infty}{(pq\kappa^{\pm1}w/z_a;p)_\infty}\nn\\
&&\times\left(x^+_j(z_1)x^+_j(z_2)x^+_{j\pm1}(w)
-(q+q^{-1})\frac{(pq^{-1}\kappa^{\mp1}z_2/w;p)_\infty}{(pq\kappa^{\mp1}z_2/w;p)_\infty}
\frac{(pq\kappa^{\pm1}w/z_2;p)_\infty}{(pq^{-1}\kappa^{\pm1}w/z_2;p)_\infty}x^+_j(z_1)x^+_{j\pm1}(w)x^+_j(z_2)\right.\nn
\ena
\bea
&&\qquad\qquad \left.+
\prod_{a=1,2}\frac{(pq^{-1}\kappa^{\mp1}z_a/w;p)_\infty}{(pq\kappa^{\mp1}z_a/w;p)_\infty}
\frac{(pq\kappa^{\pm1}w/z_a;p)_\infty}{(pq^{-1}\kappa^{\pm1}w/z_a;p)_\infty}x^+_{j\pm1}(w)x^+_j(z_1)x^+_j(z_2)\right)\nn\\
&&=-(z_1\leftrightarrow z_2).\lb{Serrexp2}
\ena
%
We prove the case containing  $x_{j+1}(w)$. The other case is similar. 
From \eqref{ej0} we have
\bea
&&x^+_j(z_1)x^+_j(z_2)x^+_{j+1}(w)\ket{\la}^{(k)}_u/C_+^{3}\nn\\
&&=q^{-1}\frac{\theta(q^2z_1/z_2)}{\theta(z_1/z_2)}
\left(\sum_{i_1,i_2,i_3\in I  \atop{\la_{i_1}+j+1\equiv i_1+k \atop {\la_{i_2}+j+2\equiv i_2+k \atop i_2=i_3<i_1}}}
\delta(u_{i_1}/z_1)\delta(q_1u_{i_2}/z_2)\delta(u_{i_3}/w)A^+_{i_2=i_3<i_1}e^{-2Q_j-Q_{j+1}}\ket{\la+\bold{1}_{i_1}+2\bold{1}_{i_2}}^{(k)}_u\right.\nn\\
&&\quad\left. +q^{-1}\frac{\theta(q^2z_2/z_1)}{\theta(z_2/z_1)}
\sum_{i_1,i_2,i_3\in I  \atop{\la_{i_1}+j+1\equiv i_1+k \atop {\la_{i_2}+j+2\equiv i_2+k \atop i_1<i_2=i_3}}}
\delta(u_{i_1}/z_1)\delta(q_1u_{i_2}/z_2)\delta(u_{i_3}/w)A^+_{i_1<i_2=i_3}e^{-2Q_j-Q_{j+1}}\ket{\la+\bold{1}_{i_1}+2\bold{1}_{i_2}}^{(k)}_u
\right)
\nn\\
&&+q\frac{\theta(q_3z_2/w)}{\theta(q_1^{-1}z_2/w)}
\left(
\sum_{i_1,i_2,i_3\in I  \atop{\la_{i_1}+j+2\equiv i_1+k \atop {\la_{i_2}+j+1\equiv i_2+k \atop i_1=i_3<i_2}}}
\delta(q_1u_{i_1}/z_1)\delta(u_{i_2}/z_2)\delta(u_{i_3}/w)A^+_{i_1=i_3<i_2}e^{-2Q_j-Q_{j+1}}\ket{\la+2\bold{1}_{i_1}+\bold{1}_{i_2}}^{(k)}_u\right.
\nn\\
&&\quad \left. +q^{-1}\frac{\theta(q^2z_1/z_2)}{\theta(z_1/z_2)}
\sum_{i_1,i_2,i_3\in I  \atop{\la_{i_1}+j+2\equiv i_1+k \atop {\la_{i_2}+j+1\equiv i_2+k \atop i_2<i_1=i_3}}}
\delta(q_1u_{i_1}/z_1)\delta(u_{i_2}/z_2)\delta(u_{i_3}/w)A^+_{i_2<i_1=i_3}e^{-2Q_j-Q_{j+1}}\ket{\la+2\bold{1}_{i_1}+\bold{1}_{i_2}}^{(k)}_u
\right)
\nn\\
&&+q\frac{\theta(q^2z_1/z_2)}{\theta(z_1/z_2)}\frac{\theta(q_3z_1/w)}{\theta(q_1^{-1}z_1/w)}\frac{\theta(q_3z_2/w)}{\theta(q_1^{-1}z_2/w)}\nn\\
&&\times \sum_{\sigma\in \gS_3}\mathop{{\sum}^*}_{
\atop  i_{\sigma(1)}<i_{\si(2)}<i_{\si(3)}
}
\delta(u_{i_1}/z_1)\delta(u_{i_2}/z_2)\delta(u_{i_3}/w)A^+_{i_{\si(1)}<i_{\si(2)}<i_{\si(3)}}
e^{-2Q_j-Q_{j+1}}\ket{\la+\bold{1}_{i_1}+\bold{1}_{i_1}+\bold{1}_{i_3}}^{(k)}_u,\nn\\
&&\lb{ejejejpo}
\ena 
where $\sum^*$ denotes a summation over $i_1,i_2,i_3\in I $ satisfying $\la_{i_1}+j+1\equiv i_1+k$,  
 $\la_{i_2}+j+1\equiv i_2+k$,  $\la_{i_3}+j+2\equiv i_3+k$, and we set
\be
 A^+_{i_2=i_3<i_1}&=&   \prod_{s_1=1\atop { \not=i_2 \atop \la_{s_1}+j\equiv s_1+k}}^{i_1-1} q^{-1}\frac{\theta(q_3^{-1}u_{i_1}/u_{s_1})}{\theta(q_1u_{i_1}/u_{s_1})}\prod_{s_1=1\atop { 
 \la_{s_1}+j+1\equiv s_1+k}}^{i_1-1} q\frac{\theta(q_1q_3u_{i_1}/u_{s_1})}{\theta(u_{i_1}/u_{s_1})}\\
&&\qquad \times\prod_{s_2=1\atop \la_{s_2}+j\equiv s_2+k}^{i_2-1} q^{-1}\frac{\theta(q_1q_3^{-1}u_{i_2}/{u}_{s_2})}{\theta(q_1^2u_{i_2}/{u}_{s_2})}\prod_{s_2=1\atop \la_{s_2}+j+1\equiv s_2+k}^{i_2-1} q\frac{\theta(q_1^2q_3u_{i_2}/{u}_{s_2})}{\theta(q_1u_{i_2}/{u}_{s_2})}\\
&&\qquad \times\prod_{s_3=1\atop { \la_{s_3}+j+1\equiv s_3+k}}^{i_2-1} q^{-1}\frac{\theta(q_3^{-1}u_{i_2}/{u}_{s_3})}{\theta(q_1u_{i_2}/{u}_{s_3})}\prod_{s_3=1\atop \la_{s_3}+j+2\equiv s_3+k}^{i_2-1} q\frac{\theta(q_1q_3u_{i_2}/{u}_{s_3})}{\theta(u_{i_2}/{u}_{s_3})},\\
A^+_{i_1<i_2=i_3}&=&\prod_{s_1=1\atop {  \la_{s_1}+j\equiv s_1+k}}^{i_1-1} q^{-1}\frac{\theta(q_3^{-1}u_{i_1}/u_{s_1})}{\theta(q_1u_{i_1}/u_{s_1})}\prod_{s_1=1\atop { 
 \la_{s_1}+j+1\equiv s_1+k}}^{i_1-1} q\frac{\theta(q_1q_3u_{i_1}/u_{s_1})}{\theta(u_{i_1}/u_{s_1})}\\
&&\qquad \times\prod_{s_2=1\atop \la_{s_2}+j\equiv s_2+k}^{i_2-1} q^{-1}\frac{\theta(q_1q_3^{-1}u_{i_2}/{u}_{s_2})}{\theta(q_1^2u_{i_2}/{u}_{s_2})}\prod_{s_2=1\atop \la_{s_2}+j+1\equiv s_2+k}^{i_2-1} q\frac{\theta(q_1^2q_3u_{i_2}/{u}_{s_2})}{\theta(q_1u_{i_2}/{u}_{s_2})}\\
&&\qquad \times\prod_{s_3=1\atop {\not=i_1 \atop \la_{s_3}+j+1\equiv s_3+k}}^{i_2-1} q^{-1}\frac{\theta(q_3^{-1}u_{i_2}/{u}_{s_3})}{\theta(q_1u_{i_2}/{u}_{s_3})}\prod_{s_3=1\atop \la_{s_3}+j+2\equiv s_3+k}^{i_2-1} q\frac{\theta(q_1q_3u_{i_2}/{u}_{s_3})}{\theta(u_{i_2}/{u}_{s_3})},\\
A^+_{i_1=i_3<i_2}&=&\prod_{s_1=1\atop {  \la_{s_1}+j\equiv s_1+k}}^{i_1-1} q^{-1}\frac{\theta(q_1q_3^{-1}u_{i_1}/u_{s_1})}{\theta(q_1^2u_{i_1}/u_{s_1})}\prod_{s_1=1\atop { 
 \la_{s_1}+j+1\equiv s_1+k}}^{i_1-1} q\frac{\theta(q_1^2q_3u_{i_1}/u_{s_1})}{\theta(q_1u_{i_1}/u_{s_1})}\\
&&\qquad \times\prod_{s_2=1\atop {\not=i_1 \atop \la_{s_2}+j\equiv s_2+k}}^{i_2-1} q^{-1}\frac{\theta(q_3^{-1}u_{i_2}/{u}_{s_2})}{\theta(q_1u_{i_2}/{u}_{s_2})}\prod_{s_2=1\atop \la_{s_2}+j+1\equiv s_2+k}^{i_2-1} q\frac{\theta(q_1q_3u_{i_2}/{u}_{s_2})}{\theta(u_{i_2}/{u}_{s_2})}\\
&&\qquad \times\prod_{s_3=1\atop { \la_{s_3}+j+1\equiv s_3+k}}^{i_1-1} q^{-1}\frac{\theta(q_3^{-1}u_{i_1}/{u}_{s_3})}{\theta(q_1u_{i_1}/{u}_{s_3})}\prod_{s_3=1\atop \la_{s_3}+j+2\equiv s_3+k}^{i_1-1} q\frac{\theta(q_1q_3u_{i_1}/{u}_{s_3})}{\theta(u_{i_1}/{u}_{s_3})},\\
A^+_{i_2<i_1=i_3}&=&\prod_{s_1=1\atop {\not=i_2 \atop  \la_{s_1}+j\equiv s_1+k}}^{i_1-1} q^{-1}\frac{\theta(q_1q_3^{-1}u_{i_1}/u_{s_1})}{\theta(q_1^2u_{i_1}/u_{s_1})}\prod_{s_1=1\atop { 
 \la_{s_1}+j+1\equiv s_1+k}}^{i_1-1} q\frac{\theta(q_1^2q_3u_{i_1}/u_{s_1})}{\theta(q_1u_{i_1}/u_{s_1})}\\
&&\qquad \times\prod_{s_2=1\atop { \atop \la_{s_2}+j\equiv s_2+k}}^{i_2-1} q^{-1}\frac{\theta(q_3^{-1}u_{i_2}/{u}_{s_2})}{\theta(q_1u_{i_2}/{u}_{s_2})}\prod_{s_2=1\atop \la_{s_2}+j+1\equiv s_2+k}^{i_2-1} q\frac{\theta(q_1q_3u_{i_2}/{u}_{s_2})}{\theta(u_{i_2}/{u}_{s_2})}\\
&&\qquad \times\prod_{s_3=1\atop {\not=i_2 \atop \la_{s_3}+j+1\equiv s_3+k}}^{i_1-1} q^{-1}\frac{\theta(q_3^{-1}u_{i_1}/{u}_{s_3})}{\theta(q_1u_{i_1}/{u}_{s_3})}\prod_{s_3=1\atop \la_{s_3}+j+2\equiv s_3+k}^{i_1-1} q\frac{\theta(q_1q_3u_{i_1}/{u}_{s_3})}{\theta(u_{i_1}/{u}_{s_3})},
\en
\be
A^+_{i_{1}<i_{2}<i_{3}}&=&\prod_{s_1=1\atop {
 \la_{s_1}+j\equiv s_1+k}}^{i_1-1} q^{-1}\frac{\theta(q_3^{-1}u_{i_1}/u_{s_1})}{\theta(q_1u_{i_1}/u_{s_1})}\prod_{s_1=1\atop { 
 \la_{s_1}+j+1\equiv s_1+k}}^{i_1-1} q\frac{\theta(q_1q_3u_{i_1}/u_{s_1})}{\theta(u_{i_1}/u_{s_1})}\\
&&\qquad \times\prod_{s_2=1\atop { \atop \la_{s_2}+j\equiv s_2+k}}^{i_2-1} q^{-1}\frac{\theta(q_3^{-1}u_{i_2}/{u}_{s_2})}{\theta(q_1u_{i_2}/{u}_{s_2})}\prod_{s_2=1\atop 
{\not=i_1 \atop \la_{s_2}+j+1\equiv s_2+k}
}^{i_2-1} q\frac{\theta(q_1q_3u_{i_2}/{u}_{s_2})}{\theta(u_{i_2}/{u}_{s_2})}\\
&&\qquad \times\prod_{s_3=1\atop 
{\not=i_1, i_2 \atop \la_{s_3}+j+1\equiv s_3+k}}^{i_3-1} q^{-1}\frac{\theta(q_3^{-1}u_{i_3}/{u}_{s_3})}{\theta(q_1u_{i_3}/{u}_{s_3})}\prod_{s_3=1\atop \la_{s_3}+j+2\equiv s_3+k}^{i_3-1} q\frac{\theta(q_1q_3u_{i_3}/{u}_{s_3})}{\theta(u_{i_3}/{u}_{s_3})},
\\
A^+_{i_{1}<i_{3}<i_{2}}&=&\prod_{s_1=1\atop {
 \la_{s_1}+j\equiv s_1+k}}^{i_1-1} q^{-1}\frac{\theta(q_3^{-1}u_{i_1}/u_{s_1})}{\theta(q_1u_{i_1}/u_{s_1})}\prod_{s_1=1\atop { 
 \la_{s_1}+j+1\equiv s_1+k}}^{i_1-1} q\frac{\theta(q_1q_3u_{i_1}/u_{s_1})}{\theta(u_{i_1}/u_{s_1})}\\
&&\qquad \times\prod_{s_2=1\atop { \atop \la_{s_2}+j\equiv s_2+k}}^{i_2-1} q^{-1}\frac{\theta(q_3^{-1}u_{i_2}/{u}_{s_2})}{\theta(q_1u_{i_2}/{u}_{s_2})}\prod_{s_2=1\atop 
{\not=i_1 \atop \la_{s_2}+j+1\equiv s_2+k}
}^{i_2-1} q\frac{\theta(q_1q_3u_{i_2}/{u}_{s_2})}{\theta(u_{i_2}/{u}_{s_2})}\\
&&\qquad \times\prod_{s_3=1\atop 
{\not=i_1
 \atop \la_{s_3}+j+1\equiv s_3+k}}^{i_3-1} q^{-1}\frac{\theta(q_3^{-1}u_{i_3}/{u}_{s_3})}{\theta(q_1u_{i_3}/{u}_{s_3})}\prod_{s_3=1\atop \la_{s_3}+j+2\equiv s_3+k}^{i_3-1} q\frac{\theta(q_1q_3u_{i_3}/{u}_{s_3})}{\theta(u_{i_3}/{u}_{s_3})},
\\
A^+_{i_{3}<i_{1}<i_{2}}&=&\prod_{s_1=1\atop {
 \la_{s_1}+j\equiv s_1+k}}^{i_1-1} q^{-1}\frac{\theta(q_3^{-1}u_{i_1}/u_{s_1})}{\theta(q_1u_{i_1}/u_{s_1})}\prod_{s_1=1\atop { 
 \la_{s_1}+j+1\equiv s_1+k}}^{i_1-1} q\frac{\theta(q_1q_3u_{i_1}/u_{s_1})}{\theta(u_{i_1}/u_{s_1})}\\
&&\qquad \times\prod_{s_2=1\atop { \atop \la_{s_2}+j\equiv s_2+k}}^{i_2-1} q^{-1}\frac{\theta(q_3^{-1}u_{i_2}/{u}_{s_2})}{\theta(q_1u_{i_2}/{u}_{s_2})}\prod_{s_2=1\atop 
{\not=i_1 \atop \la_{s_2}+j+1\equiv s_2+k}
}^{i_2-1} q\frac{\theta(q_1q_3u_{i_2}/{u}_{s_2})}{\theta(u_{i_2}/{u}_{s_2})}\\
&&\qquad \times\prod_{s_3=1\atop 
{
 \atop \la_{s_3}+j+1\equiv s_3+k}}^{i_3-1} q^{-1}\frac{\theta(q_3^{-1}u_{i_3}/{u}_{s_3})}{\theta(q_1u_{i_3}/{u}_{s_3})}\prod_{s_3=1\atop \la_{s_3}+j+2\equiv s_3+k}^{i_3-1} q\frac{\theta(q_1q_3u_{i_3}/{u}_{s_3})}{\theta(u_{i_3}/{u}_{s_3})},
\\
A^+_{i_{2}<i_{1}<i_{3}}&=&\prod_{s_1=1\atop {
\not=i_2 \atop 
 \la_{s_1}+j\equiv s_1+k}}^{i_1-1} q^{-1}\frac{\theta(q_3^{-1}u_{i_1}/u_{s_1})}{\theta(q_1u_{i_1}/u_{s_1})}\prod_{s_1=1\atop { 
 \la_{s_1}+j+1\equiv s_1+k}}^{i_1-1} q\frac{\theta(q_1q_3u_{i_1}/u_{s_1})}{\theta(u_{i_1}/u_{s_1})}\\
&&\qquad \times\prod_{s_2=1\atop { \atop \la_{s_2}+j\equiv s_2+k}}^{i_2-1} q^{-1}\frac{\theta(q_3^{-1}u_{i_2}/{u}_{s_2})}{\theta(q_1u_{i_2}/{u}_{s_2})}\prod_{s_2=1\atop 
{
 \la_{s_2}+j+1\equiv s_2+k}
}^{i_2-1} q\frac{\theta(q_1q_3u_{i_2}/{u}_{s_2})}{\theta(u_{i_2}/{u}_{s_2})}\\
&&\qquad \times\prod_{s_3=1\atop 
{\not=i_1, i_2 \atop \la_{s_3}+j+1\equiv s_3+k}}^{i_3-1} q^{-1}\frac{\theta(q_3^{-1}u_{i_3}/{u}_{s_3})}{\theta(q_1u_{i_3}/{u}_{s_3})}\prod_{s_3=1\atop \la_{s_3}+j+2\equiv s_3+k}^{i_3-1} q\frac{\theta(q_1q_3u_{i_3}/{u}_{s_3})}{\theta(u_{i_3}/{u}_{s_3})},
\en
\be
A^+_{i_{2}<i_{3}<i_{1}}&=&\prod_{s_1=1\atop {
\not=i_2 \atop 
 \la_{s_1}+j\equiv s_1+k}}^{i_1-1} q^{-1}\frac{\theta(q_3^{-1}u_{i_1}/u_{s_1})}{\theta(q_1u_{i_1}/u_{s_1})}\prod_{s_1=1\atop { 
 \la_{s_1}+j+1\equiv s_1+k}}^{i_1-1} q\frac{\theta(q_1q_3u_{i_1}/u_{s_1})}{\theta(u_{i_1}/u_{s_1})}\\
&&\qquad \times\prod_{s_2=1\atop { \atop \la_{s_2}+j\equiv s_2+k}}^{i_2-1} q^{-1}\frac{\theta(q_3^{-1}u_{i_2}/{u}_{s_2})}{\theta(q_1u_{i_2}/{u}_{s_2})}\prod_{s_2=1\atop 
{\not=i_1 \atop \la_{s_2}+j+1\equiv s_2+k}
}^{i_2-1} q\frac{\theta(q_1q_3u_{i_2}/{u}_{s_2})}{\theta(u_{i_2}/{u}_{s_2})}\\
&&\qquad \times\prod_{s_3=1\atop 
{\not=i_2
 \atop \la_{s_3}+j+1\equiv s_3+k}}^{i_3-1} q^{-1}\frac{\theta(q_3^{-1}u_{i_3}/{u}_{s_3})}{\theta(q_1u_{i_3}/{u}_{s_3})}\prod_{s_3=1\atop \la_{s_3}+j+2\equiv s_3+k}^{i_3-1} q\frac{\theta(q_1q_3u_{i_3}/{u}_{s_3})}{\theta(u_{i_3}/{u}_{s_3})},
\\
A^+_{i_{3}<i_{2}<i_{1}}&=&\prod_{s_1=1\atop {
\not=i_2 \atop 
 \la_{s_1}+j\equiv s_1+k}}^{i_1-1} q^{-1}\frac{\theta(q_3^{-1}u_{i_1}/u_{s_1})}{\theta(q_1u_{i_1}/u_{s_1})}\prod_{s_1=1\atop { 
 \la_{s_1}+j+1\equiv s_1+k}}^{i_1-1} q\frac{\theta(q_1q_3u_{i_1}/u_{s_1})}{\theta(u_{i_1}/u_{s_1})}\\
&&\qquad \times\prod_{s_2=1\atop { \atop \la_{s_2}+j\equiv s_2+k}}^{i_2-1} q^{-1}\frac{\theta(q_3^{-1}u_{i_2}/{u}_{s_2})}{\theta(q_1u_{i_2}/{u}_{s_2})}\prod_{s_2=1\atop 
{\not=i_1 \atop \la_{s_2}+j+1\equiv s_2+k}
}^{i_2-1} q\frac{\theta(q_1q_3u_{i_2}/{u}_{s_2})}{\theta(u_{i_2}/{u}_{s_2})}\\
&&\qquad \times\prod_{s_3=1\atop 
{
 \atop \la_{s_3}+j+1\equiv s_3+k}}^{i_3-1} q^{-1}\frac{\theta(q_3^{-1}u_{i_3}/{u}_{s_3})}{\theta(q_1u_{i_3}/{u}_{s_3})}\prod_{s_3=1\atop \la_{s_3}+j+2\equiv s_3+k}^{i_3-1} q\frac{\theta(q_1q_3u_{i_3}/{u}_{s_3})}{\theta(u_{i_3}/{u}_{s_3})}.
\en 
Similarly we obtain 
\bea
&&x^+_j(z_1)x^+_{j+1}(w)x^+_j(z_2)\ket{\la}^{(k)}_u/C_+^{3}\nn\\
&&=
\sum_{i,i',i''\in I  \atop{\la_i+j+1\equiv i+k \atop {\la_\ip+j+2\equiv ip+k \atop i'=i''<i}}}
\delta(q_1u_\ip/z_1)\delta(u_i/z_2)\delta(u_\ipp/w)A^+_{i'=i''<i}e^{-2Q_j-Q_{j+1}}\ket{\la+\bold{1}_i+2\bold{1}_\ip}^{(k)}_u\nn\\
&&\qquad+q^{-1}\frac{\theta(q^2z_1/z_2)}{\theta(z_1/z_2)}\sum_{i,i',i''\in I  \atop{\la_i+j+1\equiv i+k \atop {\la_\ip+j+2\equiv ip+k \atop i<i'=i''}}}
\delta(q_1u_\ip/z_1)\delta(u_i/z_2)\delta(u_\ipp/w)A^+_{i<i'=i''}e^{-2Q_j-Q_{j+1}}\ket{\la+\bold{1}_i+2\bold{1}_\ip}^{(k)}_u\nn
\\
&&\qquad+\frac{\theta(q^2z_1/z_2)}{\theta(z_1/z_2)}\frac{\theta(q_3z_1/w)}{\theta(q_1^{-1}z_1/w)}\nn\\
&&\qquad \times \sum_{\sigma\in \gS_3}\mathop{{\sum}^*}_{
\atop  i_{\sigma(1)}<i_{\si(2)}<i_{\si(3)}
}
\delta(u_{i_1}/z_1)\delta(u_{i_2}/z_2)\delta(u_{i_3}/w)A^+_{i_{\si(1)}<i_{\si(2)}<i_{\si(3)}}
e^{-2Q_j-Q_{j+1}}\ket{\la+\bold{1}_{i_1}+\bold{1}_{i_1}+\bold{1}_{i_3}}^{(k)}_u,\nn\\
&&\lb{ejejpoej}\\
&&x^+_{j+1}(w)x^+_j(z_1)x^+_j(z_2)\ket{\la}^{(k)}_u/C_+^{3}\nn\\
&&=
q^{-1}\frac{\theta(q^2z_1/z_2)}{\theta(z_1/z_2)}\sum_{\sigma\in \gS_3}\mathop{{\sum}^*}_{\atop  i_{\sigma(1)}<i_{\si(2)}<i_{\si(3)}
}
\delta(u_{i_1}/z_1)\delta(u_{i_2}/z_2)\delta(u_{i_3}/w)
\nn\\
&&\qquad\qquad\qquad\qquad\qquad\qquad\qquad\times A^+_{i_{\si(1)}<i_{\si(2)}<i_{\si(3)}}
e^{-2Q_j-Q_{j+1}}\ket{\la+\bold{1}_{i_1}+\bold{1}_{i_1}+\bold{1}_{i_3}}^{(k)}_u.
\lb{ejpoejej}
\ena  
Combining these, the action of the LHS of \eqref{Serrexp2} on $\ket{\la}^{(k)}_u/C^3_{+}$ is given by
\be
&&q^{-1}\frac{(pq^2z_1/z_2;p)_\infty(pq^2z_2/z_1;p)_\infty(p;p)_\infty}{z_2\theta(z_1/z_2)}
\prod_{a=1,2}\frac{(pq^{-1}\kappa w/z_a;p)_\infty}{(pq\kappa w/z_a;p)_\infty}\\
&&\times\left[
z_2(1-q^2z_1/z_2)\left(
\sum_{i_1,i_2,i_3\in I  \atop{\la_{i_1}+j+1\equiv i_1+k \atop {\la_{i_2}+j+2\equiv i_2+k \atop i_2=i_3<i_1}}}
\delta(u_{i_1}/z_1)\delta(q_1u_{i_2}/z_2)\delta(u_{i_3}/w)A^+_{i_2=i_3<i_1}e^{-2Q_j-Q_{j+1}}\ket{\la+\bold{1}_{i_1}+2\bold{1}_{i_2}}^{(k)}_u\right.\right.\\
&& \left.+q^{-1}\frac{\theta(q^2z_2/z_1)}{\theta(z_2/z_1)}
\sum_{i_1,i_2,i_3\in I  \atop{\la_{i_1}+j+1\equiv i_1+k \atop {\la_{i_2}+j+2\equiv i_2+k \atop i_1<i_2=i_3}}}
\delta(u_{i_1}/z_1)\delta(q_1u_{i_2}/z_2)\delta(u_{i_3}/w)A^+_{i_1<i_2=i_3}e^{-2Q_j-Q_{j+1}}\ket{\la+\bold{1}_{i_1}+2\bold{1}_{i_2}}^{(k)}_u\right)
\\
&&+z_1(1-q^2z_2/z_1)\left(
\sum_{i_1,i_2,i_3\in I  \atop{\la_{i_1}+j+2\equiv i_1+k \atop {\la_{i_2}+j+1\equiv i_2+k \atop i_1=i_3<i_2}}}
\delta(q_1u_{i_1}/z_1)\delta(u_{i_2}/z_2)\delta(u_{i_3}/w)A^+_{i_1=i_3<i_2}e^{-2Q_j-Q_{j+1}}\ket{\la+2\bold{1}_{i_1}+\bold{1}_{i_2}}^{(k)}_u\right.\\
&&\quad \left. +q^{-1}\frac{\theta(q^2z_1/z_2)}{\theta(z_1/z_2)}
\sum_{i_1,i_2,i_3\in I  \atop{\la_{i_1}+j+2\equiv i_1+k \atop {\la_{i_2}+j+1\equiv i_2+k \atop i_2<i_1=i_3}}}
\delta(q_1u_{i_1}/z_1)\delta(u_{i_2}/z_2)\delta(u_{i_3}/w)A^+_{i_2<i_1=i_3}e^{-2Q_j-Q_{j+1}}\ket{\la+2\bold{1}_{i_1}+\bold{1}_{i_2}}^{(k)}_u
\right)
\\
&&+q_3(1-q^2)z_1z_2\frac{\theta(q_3z_1/w)}{\theta(q_1^{-1}z_1/w)}\frac{\theta(q_3z_2/w)}{\theta(q_1^{-1}z_2/w)}
\frac{(1-q^2z_1/z_2)(1-q^2z_2/z_1)}{(1-q_3z_1)(1-q_3z_2)}\\
&&\times 
\left.
\sum_{\sigma\in \gS_3}\mathop{{\sum}^*}_{
\atop{  i_{\sigma(1)}<i_{\si(2)}<i_{\si(3)}
\atop{ \quad \atop{ \quad \atop{ \quad \atop\quad}}} }
}
\delta(u_{i_1}/z_1)\delta(u_{i_2}/z_2)\delta(u_{i_3}/w)A^+_{i_{\si(1)}<i_{\si(2)}<i_{\si(3)}}
e^{-2Q_j-Q_{j+1}}\ket{\la+\bold{1}_{i_1}+\bold{1}_{i_1}+\bold{1}_{i_3}}^{(k)}_u\right].
\en
Due to $z_2\theta(z_1/z_2)=-z_1\theta(z_2/z_1)$,  this expression gives the RHS of \eqref{Serrexp2} 
on $\ket{\la}^{(k)}_u/C_+^3$  by exchanging $z_1$ and $z_2$. 
 \qed 
  
 \end{appendix}

\newpage

\renewcommand{\baselinestretch}{0.7}

\end{document}